\documentclass[12pt]{amsart}
\usepackage{amsmath,amssymb,amscd,amsthm,array,mathrsfs}
\usepackage[a4paper,margin=23mm]{geometry}
\usepackage[hidelinks]{hyperref}
\newcommand{\A}{\mathcal{A}}
\newcommand{\Ll}{\mathrm{L}}
\newcommand{\Rr}{\mathrm{R}}
\newcommand{\Al}{\mathrm{Al}}
\newcommand{\prs}{\langle\cdot,\cdot\rangle}
\newcommand{\la}{\lambda}
\newtheorem{Theorem}{Theorem}[section]
\newtheorem{Proposition}[Theorem]{Proposition}
\newtheorem{Lemma}[Theorem]{Lemma}
\newtheorem{Corollary}[Theorem]{Corollary}
\theoremstyle{definition}
\newtheorem{Definition}[Theorem]{Definition}
\theoremstyle{remark}
\newtheorem{Remark}[Theorem]{Remark}
\newtheorem*{UnnumberedRemark}{Remark}
\allowdisplaybreaks[2]
\begin{document}
\title[Euclidean left-alternative algebras]{A Jordan-theoretic approach to the complete classification of Euclidean left-alternative algebras}

\author{Sa\"id Benayadi} 
 \address{Universit\'e de Lorraine, Laboratoire IECL, CNRS-UMR 7502, UFR MIM, 3 rue Augustin Fresnel, BP 45112, 57073 Metz Cedex 03, France}
 \email{said.benayadi@univ-lorraine.fr}

\author{Mohamed Boucetta}
\address{Universit\'e Cadi-Ayyad,
	Facult\'e des sciences et techniques,
	B.P. 549, Marrakech, Maroc.} \email{m.boucetta@uca.ac.ma}
\author{Hamza El Ouali}
\address{Division of Science and Mathematics, New York University Abu Dhabi, P.O. Box 129188, Abu Dhabi, United Arab Emirates.}
\email{hamza.el.ouali@nyu.edu}
\author{Hicham Lebzioui}
\address{Universit\'e Sultan Moulay Slimane,
	\'{E}cole Sup\'{e}rieure de Technologie-Kh\'{e}nifra, B.P : 170, Kh\'{e}nifra, Maroc. }
    \email{h.lebzioui@usms.ma}

 \thanks{H. E. O. was supported by grant NYUAD-065.}
       
       

\begin{abstract}
A pseudo-Euclidean left-alternative algebra is a left-alternative algebra endowed with a non-degenerate symmetric bilinear form with respect to which all left multiplication operators are symmetric. In this paper, we investigate finite-dimensional Euclidean left-alternative algebras, namely those for which the bilinear form is positive definite. Using the structure theory of Jordan algebras, we establish a complete structural characterization and classification of these algebras up to isometric isomorphism. In particular, we prove that every Euclidean left-alternative algebra is associative. We also investigate the Lie algebras arising from their commutator products and prove that the corresponding simply connected Riemannian Lie groups are isometric to products of a Euclidean space and real hyperbolic spaces. In particular, the associated metrics are flat if and only if the original algebras are commutative.
\end{abstract}

\maketitle

\thispagestyle{empty}
\setcounter{tocdepth}{2}
\section{Introduction}

In this paper, we consider finite-dimensional algebras over the field $\mathbb{R}$. An algebra is called \textit{trivial} if its multiplication is identically zero.

To study the curvature of a Lie group $G$ equipped with a left-invariant pseudo-Riemannian metric $\mu$, it is sufficient to reduce the study to its Lie algebra $\mathfrak{g}$, endowed with the non-degenerate symmetric bilinear form $\prs = \mu_e$, where $e$ denotes the identity element of $G$. The Levi-Civita connection of the pseudo-Riemannian manifold $(G,\mu)$ induces a bilinear product on $(\mathfrak{g},\prs)$, called the Levi-Civita product, which is given by the Koszul formula:
\begin{equation}\label{LV}
2\langle u \star v, w \rangle = \langle [u,v], w \rangle - \langle [v,w], u \rangle + \langle [w,u], v \rangle, \quad \text{for all } u,v,w \in \mathfrak{g}.
\end{equation}
This product is uniquely determined by the metric and satisfies:
\[
u \star v - v\star  u = [u,v] \quad \text{and} \quad \langle u\star v, w \rangle + \langle v, u\star w \rangle = 0.
\]
The curvature tensor of $(G, \mu)$ is given by:
\[
\mathfrak{R}(u,v) = \ell_{[u,v]} - [\ell_u, \ell_v],
\]
where $\ell_u$ denotes the left multiplication by $u$ with respect to the Levi-Civita product. The pair $(\mathfrak{g},\prs)$ is said to be flat if $\mathfrak{R} = 0$, which is equivalent to the Levi-Civita product being left-symmetric. More precisely,
\[
\mathrm{ass}(u,v,w) = \mathrm{ass}(v,u,w), \quad \text{where } \mathrm{ass}(u,v,w) = (u\star v)\star w - u\star (v\star w).
\]
Thus, flat pseudo-Riemannian Lie algebras can be regarded as left-symmetric algebras endowed with a non-degenerate symmetric bilinear form for which the left multiplications are skew-symmetric.

In \cite{Milnor}, J. Milnor characterized the Lie groups that admit flat left-invariant Riemannian metrics by showing that their Lie algebras decompose orthogonally as $\mathfrak{g} = \mathfrak{b} \oplus \mathfrak{u}$, where $\mathfrak{b}$ is an abelian subalgebra, $\mathfrak{u}$ is an abelian ideal, and the adjoint operators $\mathrm{ad}_b$ are skew-symmetric with respect to the inner product, for all $b \in \mathfrak{b}$.

Furthermore, E. I. Zelmanov in \cite{Zelmanov} studied Novikov algebras endowed with a positive-definite symmetric bilinear form such that the right multiplications are symmetric. He classified all such Novikov algebras and showed that they are necessarily associative.

Analogously, in \cite{Benayadi1} and \cite{Yan and Chen}, a generalization of the Levi-Civita product was introduced in the setting of commutative algebras with a non-degenerate symmetric bilinear form. Given a commutative algebra $(\mathcal{A}, \cdot)$ endowed with such a form $\prs$, the product $\bullet$ is defined by:
\[
2\langle u \bullet v, w \rangle = \langle u \cdot v, w \rangle - \langle v \cdot w, u \rangle + \langle w \cdot u, v \rangle,  \quad \text{for all } u,v,w \in \mathcal{A}.
\]
This product satisfies:
\[
u \bullet v + v \bullet u = u \cdot v \quad \text{and} \quad \langle u \bullet v, w \rangle = \langle v, u \bullet w \rangle.
\]
Letting $\Ll_u(v) = u \bullet v$ and $\Rr_u(v) = v \bullet u$, the curvature tensor associated with the pseudo-Euclidean Jordan algebra $(\mathcal{A}, \cdot,\prs)$ is given by:
\[
\mathcal{R}(u,v) = \Ll_{u \cdot v} - \Ll_u\circ \Ll_v - \Ll_v\circ \Ll_u.
\]
Such an algebra is said to be flat if $\mathcal{R} = 0$. Note that $\mathcal{R} = 0$ is equivalent to $(\mathcal{A},\bullet)$ being a left-alternative algebra. It was proved in \cite{Yan and Chen} that a flat pseudo-Euclidean commutative algebra is necessarily a Jordan algebra. This result is well-known; in fact, it is a special case of a more general theorem proven by Albert in \cite{Albert}.

A pseudo-Euclidean left-alternative algebra is a left-alternative algebra endowed with a non-degenerate symmetric bilinear form $\prs$ such that the left multiplication operators are symmetric with respect to this form. It turns out that if the pseudo-Euclidean Jordan algebra $(\mathcal{A}, \cdot,\prs)$ is flat, then the corresponding Levi-Civita product $\bullet$ defines a left-alternative algebra. In this case, the conditions:
\[
\langle u \bullet v, w \rangle = \langle v, u \bullet w \rangle \quad \text{and} \quad u \cdot v = u \bullet v + v \bullet u\]
are satisfied, and hence $(\mathcal{A}, \bullet,\prs)$ is a pseudo-Euclidean left-alternative algebra.

In analogy with Milnor’s theorem on flat Riemannian Lie algebras and Zelmanov’s classification of Novikov algebras endowed with a positive-definite symmetric bilinear form for which the right multiplication operators are symmetric, we aim to provide a structural characterization and classification of Euclidean left-alternative algebras. We prove, in particular, that every Euclidean left-alternative algebra is necessarily associative.

As a geometric application of our classification, we describe the simply connected Riemannian Lie groups associated with the commutator Lie algebras of Euclidean left-alternative algebras. Here a Riemannian Lie group is a Lie group endowed with a left-invariant Riemannian metric. Our classification leads to an orthogonal decomposition into an abelian ideal and ideals of the type considered by Milnor in \cite[Special Example 1.7]{Milnor}. Consequently, the associated Riemannian Lie groups are isometric, as Riemannian manifolds, to products of a Euclidean space and real hyperbolic spaces. The dimensions and sectional curvatures of the hyperbolic factors are determined explicitly by the algebraic classification parameters. In particular, the associated metric is flat if and only if the original algebra is commutative.

\textbf{This paper is organized as follows.}  In Section 2, we recall some classical results on Jordan algebras and examine the relationship between the nilpotency of the left and right multiplication operators in a left-alternative algebra. We show that these operators are simultaneously nilpotent if and only if the associated Jordan algebra is nilpotent.

In Section 3, we introduce the notion of pseudo-Euclidean left-alternative algebras and focus on the Euclidean case. We prove that the left annihilator of such an algebra coincides with the radical of its associated Jordan algebra, and that the orthogonal complement of this radical forms a semisimple subalgebra. If $(\mathcal{A},\bullet)$ is a non-associative algebra, the new non-associative algebra $(\mathcal{A},\cdot)$, where $u \cdot v = u \bullet v + v \bullet u$ for any $u,v \in \mathcal{A}$, will be denoted by $\mathcal{A}^+$. We show that if the associated Jordan algebra $\mathcal{A}^+ = (\mathcal{A}, \cdot)$ is semisimple, then the left-alternative algebra admits a unit element. In the case where $\mathcal{A}^+$ is simple, we establish that the left-alternative algebra is actually alternative.

In Section 4, we focus on the Euclidean case. We show that if a Euclidean left-alternative algebra has a simple associated Jordan algebra $\mathcal{A}^+ = (\mathcal{A}, \cdot)$, then it must be one-dimensional. We then give a full characterization and classification of Euclidean left-alternative algebras (see Theorem~\ref{Theo Milnor} and Corollary~\ref{co1}). In particular, we prove that any such algebra is associative (see Corollary~\ref{associative}). Finally, we describe the geometry of the associated simply connected Riemannian Lie groups in Theorem~\ref{hyperbolic-product}.

\section{Left-alternative algebras and Jordan algebras}

In this section, we recall some basic definitions and fundamental theorems from the theory of Jordan algebras. We then show that the left multiplication operator is nilpotent if and only if the right multiplication operator is nilpotent, for each element. We also show that all left multiplication operators are nilpotent if and only if the associated Jordan algebra is nilpotent.

Let $(\mathcal{A},\bullet)$ be a non-associative algebra. For any $u \in \mathcal{A}$, we denote by $\Ll_u$ and $\Rr_u$ the left and right multiplication operators defined respectively by:
\[
\Ll_u(v) = u \bullet v \quad \text{and} \quad \Rr_u(v) = v \bullet u, \quad \text{for all } v \in \mathcal{A}.
\]
On the underlying vector space $\mathcal{A}$, we define the corresponding anticommutator and commutator by:
\[
u \cdot v := u \bullet v + v \bullet u \quad \text{and} \quad [u,v] := u \bullet v - v \bullet u, \quad \text{for all } u,v \in \mathcal{A}.
\] 
The algebra $(\mathcal{A}, \cdot)$ (resp. $(\mathcal{A}, [\,\cdot\,,\,\cdot\,])$) will be denoted by $\mathcal{A}^+$ (resp. $\mathcal{A}^-$).

If $(\mathcal{A},\bullet)$ is a commutative  algebra, both $\Ll_u$ and $\Rr_u$ will be denoted by $\operatorname{sd}_u$ for all $u \in \mathcal{A}$.

A Jordan algebra $(\mathcal{A},\cdot)$ is a commutative  algebra satisfying the Jordan identity for all $u,v \in \mathcal{A}$:
\begin{equation}\label{Identite de Jordan}
u \cdot (u^2 \cdot v) = u^2 \cdot (u \cdot v).
\end{equation}
For more details concerning Jordan algebras, one can refer, for instance, to \cite{Faraut, Jaco, Mc, Rich}.

The algebra $(\mathcal{A}, \bullet)$ is called a \textit{Jordan-admissible algebra} if $\mathcal{A}^+$ is a Jordan algebra.

We now recall some classical results on the structure theory of Jordan algebras.
 
Let $(\mathcal{A}, \cdot)$ be a Jordan algebra. The algebra $\mathcal{A}$ is said to be \textit{solvable} if there exists an integer $n \in \mathbb{N}$ such that $\mathcal{A}^{(n)} = 0$, where the sequence $(\mathcal{A}^{(k)})_{k \geq 0}$ is defined inductively by:
\[
\mathcal{A}^{(0)} = \mathcal{A}, \quad \mathcal{A}^{(k+1)} = \mathcal{A}^{(k)} \cdot \mathcal{A}^{(k)}, \quad \forall k \geq 0.
\]
Similarly, $\mathcal{A}$ is called \textit{nilpotent} if there exists an integer $m \in \mathbb{N}$ such that $\mathcal{A}^m = 0$, where:
\[
\mathcal{A}^0 = \mathcal{A}, \quad \mathcal{A}^{k+1} = \mathcal{A} \cdot \mathcal{A}^{k}, \quad \forall k \geq 0.
\]
In contrast to the Lie algebra setting, in the theory of Jordan algebras these two notions coincide: every solvable Jordan algebra is nilpotent and vice versa (see \cite{Albert1}). Moreover, any Jordan algebra admits a unique maximal nilpotent ideal, known as its \textit{radical}, and denoted by $\mathcal{N} = \operatorname{Ra}(\mathcal{A})$.
 
A Jordan algebra is said to be \textit{simple} if it has no non-trivial proper ideals and satisfies $\mathcal{A}\cdot\mathcal{A}\neq0$. It is called \textit{semisimple} if it is isomorphic to a finite direct sum of simple Jordan algebras. In order to characterize semisimplicity, one introduces a symmetric bilinear form that plays a role analogous to the Killing form in Lie algebras. This form, known as the Albert form, is defined as follows.

\begin{Theorem}[Albert]\label{Albert-form}
Let $(\mathcal{A}, \cdot)$ be a Jordan algebra. The bilinear form $\Al : \mathcal{A} \times \mathcal{A} \to \mathbb{R}$ defined by:
\begin{equation}\label{Associative de Jordan 1}
\Al(u, v) = \operatorname{Tr}\left(\operatorname{sd}_{u \cdot v}\right)
\end{equation}
is associative in the sense that:
\[
\Al(u \cdot v, w) = \Al(u, v \cdot w), \quad \text{ for all } u, v, w \in \mathcal{A}.
\]
Moreover, the radical of the algebra coincides with the radical of this form, i.e.,
\[
\operatorname{Ra}(\mathcal{A}) = \{u \in \mathcal{A} \mid \Al(u, v) = 0, \ \forall v \in \mathcal{A} \}.
\]
\end{Theorem}
 
The fundamental structural result for Jordan algebras is the following decomposition theorem:

\begin{Theorem}[Wedderburn Principal Theorem]\label{theo decom}
Let $(\mathcal{A},\cdot)$ be a real Jordan algebra. Then there exists a direct sum decomposition of vector spaces:
\[
\mathcal{A} = \mathcal{S} \oplus \mathcal{N},
\]
where $\mathcal{N} = \operatorname{Ra}(\mathcal{A})$ is the radical of $(\mathcal{A},\cdot)$, and $\mathcal{S}$ is a maximal semisimple subalgebra isomorphic to the quotient $\mathcal{A}/\mathcal{N}$.\end{Theorem}

\begin{Definition} 
Let $(\A, \bullet)$ be an algebra.  
   \begin{enumerate}
       \item $(\A, \bullet)$ is called a left-alternative  algebra if the following identity is satisfied
\begin{equation}
\operatorname{ass}(u, v, w)=-\operatorname{ass}(v, u, w)\label{Left symmetric},
\end{equation} 
 \item $\A$ is called a right-alternative algebra if the following identity is satisfied 
\begin{equation}
\operatorname{ass}(u, v, w)=-\operatorname{ass}(u, w, v)\label{Rigth symmetric},
\end{equation} 
where $\operatorname{ass}(u, v, w)=(u\bullet v) \bullet w-u\bullet(v \bullet w)$ for any $u,v,w\in\A$.
\item An algebra is alternative if it is both left-alternative and right-alternative.
 \end{enumerate}
 
\end{Definition}

\begin{UnnumberedRemark}
	\begin{enumerate}
 \item   $(\A, \bullet)$ is a left-alternative algebra (resp. right-alternative algebra) if and only if
$$\operatorname{ass}(u, u, w)=0,\; (\text{respectively, }\operatorname{ass}(u, v, v)=0)$$
for any $u,v,w\in \A.$

\item Then, the left-alternative identity~\eqref{Left symmetric} is equivalent to either of the following conditions:
\begin{align} \label{Mulgauch}
	\Ll_{u \bullet v + v \bullet u} &= \Ll_u \circ \Ll_v + \Ll_v \circ \Ll_u, \\ \label{Muladroite}
	\Rr_{u \bullet v} - \Rr_v \circ \Rr_u &= [\Rr_v, \Ll_u],
\end{align}
for all $u, v \in \A$
\end{enumerate}
\end{UnnumberedRemark}

A. A. Albert established in \cite{Albert} that every right-alternative algebra is Jordan-admissible. Furthermore, since the opposite algebra of a right-alternative algebra is a left-alternative algebra, Albert’s theorem can be equivalently reformulated as follows:

\begin{Theorem}[\cite{Albert}]\label{Albert}
Let \((\A, \bullet)\) be a left-alternative algebra. Then \((\A, \bullet)\) is Jordan-admissible; that is, \(\A^+\) is a Jordan algebra.
\end{Theorem}

\begin{Lemma}\label{Lemme1}
Let $(\mathcal{A},\bullet)$ be a left-alternative algebra. Then, for any $n \in \mathbb{N}^*$ and all $u \in \mathcal{A}$:
\begin{equation}\label{Re4}
\Rr_{\Rr^n_{u}(u)} = \Rr_u^{n+1} + \sum_{k=1}^{n} \Rr_u^{n-k} \circ \left[\Rr_u, \Ll_{\Rr_u^{k-1}(u)}\right].
\end{equation}
\end{Lemma}
  
\begin{proof}
We establish the relation \eqref{Re4} by induction on $n$. Let $u \in \mathcal{A}$. By virtue of relation \eqref{Muladroite}, the assertion holds for $n=1$. Now, assume that the induction hypothesis holds for a given $n \ge 1$. Applying relation \eqref{Muladroite} once again, we obtain:
\[
\begin{aligned}
\Rr_{\Rr^{n+1}_{u}(u)} &= \Rr_{\Rr^{n}_u(u) \bullet u} \\
&= \Rr_u \circ \Rr_{\Rr^{n}_u(u)} + \left[\Rr_u, \Ll_{\Rr^{n}_u(u)}\right] \\
&= \Rr_u \circ \left( \Rr_u^{n+1} + \Rr_u^{n-1} \circ \left[\Rr_u, \Ll_u\right] + \dots + \left[\Rr_u, \Ll_{\Rr_u^{n-1}(u)}\right] \right) + \left[\Rr_u, \Ll_{\Rr^{n}_u(u)}\right] \\
&= \Rr_u^{n+2} + \Rr_u^n \circ \left[\Rr_u, \Ll_u\right] + \Rr_u^{n-1} \circ \left[\Rr_u, \Ll_{u\bullet u}\right] + \dots + \left[\Rr_u, \Ll_{\Rr^{n}_u(u)}\right].
\end{aligned}
\]
This completes the induction step.
\end{proof}

\begin{Proposition}
Let \((\mathcal{A}, \bullet)\) be a left-alternative algebra. Then, \(\operatorname{Tr}(\Rr_u) = 0\) for all \(u \in \mathcal{A}\) if and only if \(\Rr_u\) is nilpotent for all \(u \in \mathcal{A}\).
\end{Proposition}

\begin{proof}
Assume that \(\operatorname{Tr}(\Rr_u) = 0\) for all \(u \in \mathcal{A}\). By virtue of Lemma~\ref{Lemme1}, it follows that \(\operatorname{Tr}(\Rr_u^{n+1}) = 0\) for all $n \in \mathbb{N}^*$, since \(\operatorname{Tr}(\Rr_{\Rr_u^n(u)}) = 0\) and \(\operatorname{Tr}(\Rr_u^k \circ [\Rr_u, \Ll_v]) = 0\) for all \(k \in \mathbb{N}\) and all \(v\in\A\), by cyclicity of the trace. Therefore, we deduce that \(\Rr_u\) is nilpotent for all \(u \in \A\). Conversely, if \(\Rr_u\) is nilpotent for any \(u \in \A\), then \(\operatorname{Tr}(\Rr_u) = 0\) for all \(u \in \A\).
\end{proof} 

\begin{Lemma}\label{Lemme2}
Let $(\mathcal{A},\bullet)$ be a left-alternative algebra. Then, for any $u \in \mathcal{A}$ and all $n \in \mathbb{N}^*$:
\begin{equation}\label{Re5}
\Ll^{n}_u(u) = \Rr_u(\Ll^{n-1}_u(u)) \quad \text{and} \quad \Ll^{n}_u(u) = \Rr^{n}_u(u).
\end{equation}
\end{Lemma}

\begin{proof}
We proceed by induction on \(n\). For \(n=1\), the result is immediate. For \(n=2\), using the left-alternative identity, we have:
\[
\Ll_u^2(u) = u \bullet (u \bullet u) = (u \bullet u) \bullet u = \Rr_u(\Ll_u(u)),
\]
which validates both assertions of the lemma for \(n=2\).

Now, suppose that for some \(n \geq 2\), the following induction hypotheses hold:
\[
\Ll_u^{n-1}(u) = \Rr_u(\Ll_u^{n-2}(u)) \quad \text{and} \quad \Ll_u^{n}(u) = \Rr_u(\Ll_u^{n-1}(u)).
\]
We aim to show that \(\Ll_u^{n+1}(u) = \Rr_u(\Ll_u^n(u))\). Recall from relation \eqref{Mulgauch} that \(\Ll_{u \bullet u} = \Ll_u^2\), and from \eqref{Muladroite} we have:
\[
\Rr_{(u \bullet u) \bullet u} = \Rr_u \circ \Rr_{u \bullet u} + [\Rr_u, \Ll_{u \bullet u}],
\]
expanding this identity yields:
\[
\Rr_{(u \bullet u) \bullet u} = \Rr_u^3 + \Rr_u^2 \circ \Ll_u - \Rr_u \circ \Ll_u \circ \Rr_u + \Rr_u \circ \Ll_u^2 - \Ll_u^2 \circ \Rr_u.
\]
Similarly, we obtain:
\[
\begin{aligned}
\Rr_{u \bullet (u \bullet u)} &= \Rr_{u \bullet u} \circ \Rr_u + [\Rr_{u \bullet u}, \Ll_u] \\
&= (\Rr_u^2 + [\Rr_u, \Ll_u]) \circ \Rr_u + [\Rr_u^2 + [\Rr_u, \Ll_u], \Ll_u] \\
&= \Rr_u^3 + \Rr_u \circ \Ll_u \circ \Rr_u - \Ll_u \circ \Rr_u^2 + \Rr_u^2 \circ \Ll_u\\
&\quad - \Ll_u \circ \Rr_u^2 + \Rr_u \circ \Ll_u^2 - 2 \Ll_u \circ \Rr_u \circ \Ll_u + \Ll_u^2 \circ \Rr_u.
\end{aligned}
\]

We now evaluate these operators on the element \(\Ll_u^{n-2}(u)\). First, we find:
\[
\begin{aligned}
\Rr_{(u \bullet u) \bullet u}(\Ll_u^{n-2}(u)) &= \Rr_u^3(\Ll_u^{n-2}(u)) + \Rr_u^2(\Ll_u(\Ll_u^{n-2}(u))) - \Rr_u(\Ll_u(\Rr_u(\Ll_u^{n-2}(u)))) \\
&\quad + \Rr_u(\Ll_u^2(\Ll_u^{n-2}(u))) - \Ll_u^2(\Rr_u(\Ll_u^{n-2}(u))) \\
&= \Rr_u(\Ll_u^n(u)) + \Rr_u(\Ll_u^n(u)) - \Rr_u(\Ll_u^n(u)) + \Rr_u(\Ll_u^n(u)) - \Ll_u^{n+1}(u) \\
&= 2 \Rr_u(\Ll_u^n(u)) - \Ll_u^{n+1}(u).
\end{aligned}
\]
Similarly, the second expansion gives:
\[
\begin{aligned}
\Rr_{u \bullet (u \bullet u)}(\Ll_u^{n-2}(u)) &= \Rr_u(\Ll_u^n(u)) + \Rr_u(\Ll_u^n(u)) - \Ll_u^{n+1}(u) + \Rr_u(\Ll_u^n(u)) - \Ll_u^{n+1}(u) \\
&\quad + \Rr_u(\Ll_u^n(u)) - 2 \Ll_u^{n+1}(u) + \Ll_u^{n+1}(u) \\
&= 4 \Rr_u(\Ll_u^n(u)) - 3 \Ll_u^{n+1}(u).
\end{aligned}
\]
Since \((u \bullet u) \bullet u = u \bullet (u \bullet u)\) due to the left-alternative property, we conclude that \(\Rr_{(u \bullet u) \bullet u} = \Rr_{u \bullet (u \bullet u)}\). Hence:
\[
2 \Rr_u(\Ll_u^n(u)) - \Ll_u^{n+1}(u) = 4 \Rr_u(\Ll_u^n(u)) - 3 \Ll_u^{n+1}(u),
\]
which simplifies directly to:
\[
\Rr_u(\Ll_u^n(u)) = \Ll_u^{n+1}(u).
\]
This establishes the first identity in \eqref{Re5} by induction. The second identity follows immediately by successive applications of the first:
\[
\Ll_u^n(u) = \Rr_u(\Ll_u^{n-1}(u)) = \Rr_u^2(\Ll_u^{n-2}(u)) = \dots = \Rr_u^n(u).
\]
This completes the proof.
\end{proof}

\begin{Lemma}\label{Lemme3}
Let $(\mathcal{A}, \bullet)$ be a left-alternative algebra. Then, for any $u \in \mathcal{A}$ and all $n \in \mathbb{N}^*$, we have:
\begin{equation}\label{Re6}
\Ll_{\Ll_u^{n-1}(u)} = \Ll_u^n.
\end{equation}
\end{Lemma}

\begin{proof}
We proceed by induction on $n$. For $n = 1$, the identity holds trivially. Now, assume that the identity $\Ll_{\Ll_u^{n-1}(u)} = \Ll_u^n$ holds for a given $n \geq 1$. From Lemma~\ref{Lemme2}, we have $\Ll_u^n(u) = \Rr_u(\Ll_u^{n-1}(u))$, which implies that $u \bullet \Ll_u^{n-1}(u) = \Ll_u^{n-1}(u) \bullet u$. Consequently,
\[
2\Ll_{\Ll_u^n(u)} = 2\Ll_{u \bullet \Ll_u^{n-1}(u)} = \Ll_{u \bullet \Ll_u^{n-1}(u)} + \Ll_{\Ll_u^{n-1}(u) \bullet u}.
\]
By virtue of identity \eqref{Mulgauch}, we obtain:
\[
2\Ll_{\Ll_u^n(u)} = \Ll_u \circ \Ll_{\Ll_u^{n-1}(u)} + \Ll_{\Ll_u^{n-1}(u)} \circ \Ll_u.
\]
Applying the induction hypothesis $\Ll_{\Ll_u^{n-1}(u)} = \Ll_u^n$, it follows that:
\[
2\Ll_{\Ll_u^n(u)} = \Ll_u \circ \Ll_u^n + \Ll_u^n \circ \Ll_u = 2\Ll_u^{n+1}.
\]
Thus, we get  $\Ll_{\Ll_u^n(u)} = \Ll_u^{n+1}$, which completes the induction.
\end{proof}

\begin{Proposition}\label{Nil L and R}
Let $(\mathcal{A}, \bullet)$ be a left-alternative algebra and let $u \in \mathcal{A}$. Then, $\Ll_u$ is nilpotent if and only if $\Rr_u$ is nilpotent.
\end{Proposition}

\begin{proof}
Let $u \in \mathcal{A}$. First, assume that $\Rr_u$ is nilpotent. By Lemma~\ref{Lemme2}, we have $\Ll_u^{n}(u) = \Rr_u^{n}(u)$ for all $n \geq 1$, and by Lemma~\ref{Lemme3}, we obtain $\Ll_{\Ll_{u}^{n}(u)} = \Ll_u^{n+1}$. Consequently, $\Ll_{\Rr_{u}^{n}(u)} = \Ll_u^{n+1}$ for all $n \geq 1$. Since $\Rr_u$ is nilpotent, there exists $n_0 \in \mathbb{N}^*$ such that $\Rr_u^{n_0} = 0$, which implies $\Rr_u^{n_0}(u) = 0$. Therefore, $\Ll_u^{n_0+1} =\Ll_{\Rr_{u}^{n_0}(u)} = 0$, proving that $\Ll_u$ is nilpotent.

Conversely, assume that $\Ll_u$ is nilpotent, so $\Ll_u^{n_0} = 0$ for some $n_0 \in \mathbb{N}^*$. By Lemma~\ref{Lemme2}, we have $\Rr_u^{n_0}(u) = \Ll_u^{n_0}(u) = 0$. Now, utilizing Lemma~\ref{Lemme1}, we have for all $n \geq n_0$:
\[
\Rr_{\Rr_u^n(u)} = \Rr_u^{n+1} + \Rr_u^{n-1} \circ [\Rr_u, \Ll_u] + \dots + [\Rr_u, \Ll_{\Rr_u^{n-1}(u)}].
\]
Since $\Rr_u^n(u) = 0$ for all $n \geq n_0$, the left-hand side vanishes, meaning $\Rr_{\Rr_u^n(u)} = 0$.  Moreover, since \(\operatorname{Tr}(\Rr_u^k \circ [\Rr_u, \Ll_v]) = 0\) for all \(k \in \mathbb{N}\) and all \(v \in \A\), we deduce that \(\operatorname{Tr}(\Rr_u^{n+1}) = 0\) for all \(n \geq n_0\), and hence $\Rr_u$ is nilpotent.
\end{proof}

We now obtain a useful characterization of nilpotency for Jordan algebras obtained from left-alternative algebras:

\begin{Theorem}\label{Nilpotent algebra}
Let $(\mathcal{A}, \bullet)$ be a left-alternative algebra. Then, the Jordan algebra $\mathcal{A}^+$ is nilpotent if and only if the left multiplication operator $\Ll_u$ is nilpotent for all $u \in \mathcal{A}$.
\end{Theorem}

\begin{proof}
Assume first that $\Ll_u$ is nilpotent for all $u \in \mathcal{A}$. By Proposition~\ref{Nil L and R}, the right multiplication operator $\Rr_u$ is also nilpotent for all $u$. Hence, for all $u,v \in \mathcal{A}$, we have:
\[
\Al(u,v) = \operatorname{tr}(\operatorname{sd}_{u \cdot v}) = \operatorname{tr}(\Ll_{u \cdot v}) + \operatorname{tr}(\Rr_{u \cdot v}) = 0.
\]
Since the radical $\mathcal{N}$ of a Jordan algebra coincides with the kernel of its Albert form $\Al$ (see \cite{Rich}), it follows that $\mathcal{A}^+ = \mathcal{N}$, and thus $\mathcal{A}^+$ is nilpotent.

Conversely, suppose that the Jordan algebra $\mathcal{A}^+$ is nilpotent. We show by induction that:
\begin{equation}\label{Re7}
\operatorname{sd}_u^n(u) = 2^n \Ll_u^n(u).
\end{equation}
For $n=1$, we have $\operatorname{sd}_u(u) = 2 u \bullet u = 2 \Ll_u(u)$, so the property holds. Assume it holds for some $n \geq 1$. From Lemma~\ref{Lemme2}, we have $\Ll_u^{n+1}(u) = \Rr_u(\Ll_u^n(u))$. Then:
\[
\begin{aligned}
\operatorname{sd}_u^{n+1}(u)&=\operatorname{sd}_u(\operatorname{sd}_u^n(u))=2^n\operatorname{sd}_u(\Ll_u^n(u))\\
&=2^n\big(u\bullet\Ll_u^n(u)+\Ll_u^n(u)\bullet u\big)\\
&=2^n(\Ll_u^{n+1}(u)+\Ll_u^{n+1}(u))=2^{n+1}\Ll_u^{n+1}(u).
\end{aligned}
\]
Since $\mathcal{A}^+$ is nilpotent, there exists $n_0 \in \mathbb{N}^*$ such that:
\[
\operatorname{sd}_u^{n_0}(u) = 2^{n_0} \Ll_u^{n_0}(u) = 0,
\]
which implies $\Ll_u^{n_0}(u) = 0$. By Lemma~\ref{Lemme3}, this yields $\Ll_u^{n_0+1} = 0$. Therefore, $\Ll_u$ is nilpotent for all $u \in \mathcal{A}$.
\end{proof}

\section{Euclidean Left-Alternative Algebras}

In this section, we study Euclidean left-alternative algebras and establish that the left annihilator of such an algebra coincides with the radical of its associated Jordan algebra. Moreover, the orthogonal complement of this radical is a semisimple subalgebra. We also prove that if the associated Jordan algebra $\mathcal{A}^+$ is semisimple, then the left-alternative algebra admits a unit element. Finally, in the case where $\mathcal{A}^+$ is simple, we show that the left-alternative algebra is alternative.

\begin{Definition}
Let $(\mathcal{A}, \cdot)$ be a non-associative commutative algebra endowed with a non-degenerate symmetric bilinear form $\prs$. The bilinear form is called a \textit{pseudo-Euclidean structure} on the commutative algebra $(\mathcal{A}, \cdot)$, and the triple $(\mathcal{A}, \cdot, \prs)$ defines a \textit{pseudo-Euclidean non-associative commutative algebra}. 

If, in addition, $\prs$ is positive-definite, $(\mathcal{A}, \cdot, \prs)$ will be called a \textit{Euclidean non-associative commutative algebra}.
\end{Definition}

In \cite{Benayadi1}, a construction analogous to the Levi-Civita product in pseudo-Euclidean Lie algebras was introduced for pseudo-Euclidean commutative non-associative algebras, providing a unique and natural product within this framework.

\begin{Proposition}[\cite{Benayadi1}]
Let $(\mathcal{A}, \cdot, \prs)$ be a pseudo-Euclidean commutative non-associative algebra. Then, there exists a unique bilinear product $\bullet$ on the underlying vector space of $\mathcal{A}$ such that, for all $u,v,w \in \mathcal{A}$, the following conditions hold:
\begin{align}
\label{sanstorsion}
u \cdot v &= u \bullet v + v \bullet u, \\ 
\label{compatible}
\langle u \bullet v, w \rangle &= \langle v, u \bullet w \rangle. 
\end{align}

More precisely, the product $\bullet$ is given by:
\begin{equation}\label{kozule}
\langle u \bullet v, w \rangle = \frac{1}{2} \big( \langle u \cdot v, w \rangle - \langle v \cdot w, u \rangle + \langle w \cdot u, v \rangle \big),
\end{equation}
for all $u,v,w \in \mathcal{A}$. This product $\bullet$ is called the \textit{Levi-Civita product} associated with $(\mathcal{A}, \cdot, \prs)$.
\end{Proposition}

Now, we introduce the definition of a pseudo-Euclidean left-alternative algebra.

\begin{Definition}
Let $(\mathcal{A}, \bullet)$ be a left-alternative algebra endowed with a non-degenerate symmetric bilinear form $\prs$. We say that $(\mathcal{A}, \bullet, \prs)$ is a \textit{pseudo-Euclidean left-alternative algebra} if the left multiplication operators are symmetric with respect to $\prs$, that is:
\[
\langle u \bullet v, w \rangle = \langle v, u \bullet w \rangle, \quad \text{for all } u,v,w \in \mathcal{A}.
\]

If, in addition, $\prs$ is positive-definite, $(\mathcal{A}, \bullet, \prs)$ will be called a \textit{Euclidean left-alternative algebra}.
\end{Definition}

\begin{Remark}
\begin{enumerate}
\item Let $(\mathcal{A}, \bullet, \prs)$ be a pseudo-Euclidean left-alternative algebra. According to Theorem~\ref{Albert}, $(\mathcal{A}, \bullet)$ is Jordan-admissible, and the associated Jordan product $\cdot$ is defined by:
\[
u \cdot v = u \bullet v + v \bullet u, \quad \text{ for all } u, v \in \mathcal{A}.
\]
Since the bilinear form $\prs$ satisfies:
\[
\langle u \bullet v, w \rangle = \langle v, u \bullet w \rangle,
\]
the uniqueness of the Levi-Civita product implies that $\bullet$ coincides with the Levi-Civita product associated with the pseudo-Euclidean Jordan algebra $(\mathcal{A}, \cdot, \prs)$.

Conversely, let $(\mathcal{A}, \cdot, \prs)$ be a pseudo-Euclidean Jordan algebra, and let $\bullet$ denote its Levi-Civita product. If $(\mathcal{A}, \bullet)$ is left-alternative, then $(\mathcal{A}, \bullet, \prs)$ is a pseudo-Euclidean left-alternative algebra, and the product $\bullet$ coincides with the Levi-Civita product of the pseudo-Euclidean Jordan algebra $(\mathcal{A}^+, \prs)$.

\item The authors of \cite{Yan and Chen} have adopted the geometric vocabulary of Lie theory in the setting of pseudo-Euclidean Jordan algebras. Since every pseudo-Euclidean Jordan algebra $(\mathcal{A}, \cdot, \prs)$ admits a Levi-Civita product $\bullet$, they have introduced a curvature operator $\mathcal{R}$ defined by:
\[
\mathcal{R}(u,v) := \Ll_{u \cdot v} - \Ll_u \circ \Ll_v - \Ll_v \circ \Ll_u, \quad \text{ for all } u,v \in \mathcal{A}.
\]
More explicitly,
\[
\mathcal{R}(u,v)(w) = (u \bullet v + v \bullet u) \bullet w - u \bullet (v \bullet w) - v \bullet (u \bullet w), \quad \text{ for all } u,v,w \in \mathcal{A}.
\]
If $\mathcal{R} = 0$, they have called $(\mathcal{A}, \cdot, \prs)$ a \textit{flat pseudo-Euclidean Jordan algebra}.

Let us remark that the condition $\mathcal{R} = 0$ is equivalent to the Levi-Civita product being left-alternative. Moreover, since its left multiplication operators are symmetric with respect to $\prs$, we conclude that $(\mathcal{A}, \bullet, \prs)$ is a pseudo-Euclidean left-alternative algebra in the sense of our definition.
\end{enumerate}
\end{Remark}

Let $(\mathcal{A}, \bullet, \prs)$ be a pseudo-Euclidean left-alternative algebra, and let $\mathcal{A}^+ = (\mathcal{A}, \cdot)$ denote the Jordan algebra associated with $(\mathcal{A}, \bullet)$. For any $u \in \mathcal{A}$, the left multiplication $\Ll_u$ is symmetric with respect to $\prs$, and the endomorphism $\operatorname{sd}_u: \mathcal{A} \to \mathcal{A}$, defined by $\operatorname{sd}_u(v) = u \cdot v$, satisfies the relation $\operatorname{sd}_u = \Ll_u + \Rr_u$.

Recall that $\mathcal{A} \cdot \mathcal{A} := \operatorname{span} \{ u \cdot v \mid u, v \in \mathcal{A} \}$ and the annihilator of $(\mathcal{A}, \cdot)$, defined by $N(\mathcal{A}) := \{ u \in \mathcal{A} \mid \operatorname{sd}_u = 0 \}$, are ideals of the Jordan algebra $\mathcal{A}^+$. We also define the left, right, and two-sided annihilators of $(\mathcal{A}, \bullet)$ by:
\[
N_\ell(\mathcal{A}) = \{ u \in \mathcal{A} \mid \Ll_u = 0 \}, \quad N_r(\mathcal{A}) = \{ u \in \mathcal{A} \mid \Rr_u = 0 \}, \quad \text{and} \quad N_{\ell r} = N_\ell(\mathcal{A}) \cap N_r(\mathcal{A}).
\]
Then, the following relations hold:
\begin{equation} \label{eq1}
(\mathcal{A} \cdot \mathcal{A})^\perp = \{ u \in \mathcal{A} \mid \Rr_u = -\Rr_u^* \}, \quad \text{and} \quad (\mathcal{A} \bullet \mathcal{A})^\perp = N_r(\mathcal{A}),
\end{equation}
where $\mathcal{A} \bullet \mathcal{A} = \operatorname{span} \{ u \bullet v \mid u, v \in \mathcal{A} \}$.

\begin{Proposition}\label{Euclidean}
Let $(\mathcal{A}, \bullet, \prs)$ be a pseudo-Euclidean left-alternative algebra. Then:
\begin{enumerate}
    \item For all $u, v \in N(\mathcal{A})$, we have $\Ll_u \circ \Ll_v = 0$. In particular, $\Ll_u^2 = 0$.
 \item If $\mathcal{A} \cdot \mathcal{A}$ is Euclidean (i.e., $\prs_{\vert_{(\mathcal{A} \cdot \mathcal{A}) \times (\mathcal{A} \cdot \mathcal{A})}}$ is positive-definite), then for any $u \in (\mathcal{A} \cdot \mathcal{A})^{\perp}$, we have $\Ll_u = \Rr_u = 0$.
    \item If $\prs$ is Euclidean, then $\Ll_u = \Rr_u = 0$ for all $u \in N(\mathcal{A})$.
\end{enumerate}
\end{Proposition}

\begin{proof}
\begin{enumerate}
   \item Let \(u, v \in N(\mathcal{A})\). Then, for any \(w \in \mathcal{A}\), by \eqref{Mulgauch} we have:
   \[
   \Ll_u \circ \Ll_w(v) = - \Ll_w \circ \Ll_u(v).
   \]
   Using \eqref{kozule}, we get \(u \bullet v = 0\) and \(w \bullet v = - v \bullet w\), from which it follows that:
   \[
   \Ll_u \circ \Ll_v = 0.
   \]
   In particular, taking \(u = v\) yields \(\Ll_u^2 = 0\).

   \item Let \(u \in (\mathcal{A} \cdot \mathcal{A})^{\perp}\). From \eqref{eq1}, \(\Rr_u\) is skew-symmetric with respect to $\prs$, hence:
   \[
   \langle u \bullet u, u \bullet u \rangle = -\langle u, (u \bullet u) \bullet u \rangle \overset{\eqref{Left symmetric}}{=} -\langle u, u \bullet (u \bullet u) \rangle = -\langle u \bullet u, u \bullet u \rangle,
   \]
   which implies \(\langle u \bullet u, u \bullet u \rangle = 0\). Since \(u \bullet u = \frac{1}{2} u \cdot u\), we deduce \(u \bullet u = 0\) because $\mathcal{A} \cdot \mathcal{A}$ is Euclidean. Therefore, using \eqref{kozule}, we have $2 u \bullet v = u \cdot v$ for all $v \in (\mathcal{A} \cdot \mathcal{A})^{\perp}$. Moreover, since the left multiplication operators are symmetric with respect to $\prs$, we conclude that:
   \[
   \langle u \cdot v, u \cdot v \rangle = 4 \langle u \bullet v, u \bullet v \rangle = 4 \langle v, u \bullet (u \bullet v) \rangle \overset{\eqref{Left symmetric}}{=} 4 \langle v, (u \bullet u) \bullet v \rangle = 0.
   \]
   Since \(\mathcal{A} \cdot \mathcal{A}\) is Euclidean, this implies $u \cdot v = 0$, and hence $u \bullet v =v\bullet u= 0$. It follows that:
   \[
   \Ll_u^2 = \Ll_{u \bullet u} = 0.
   \]
   By Proposition~\ref{Nil L and R}, the operator \(\Rr_u\) is also nilpotent. Set $D=\A\cdot\A$. Since $D$ is positive definite, $\A=D\oplus D^\perp$ orthogonally. We have already proved that $\Ll_u(D^\perp)=\Rr_u(D^\perp)=\{0\}$. For $x\in D$ and $y\in D^\perp$, the symmetry of $\Ll_u$ and the skew-symmetry of $\Rr_u$ give
\[
\langle\Ll_u x,y\rangle=\langle x,\Ll_u y\rangle=0,
\qquad
\langle\Rr_u x,y\rangle=-\langle x,\Rr_u y\rangle=0.
\]
Thus both operators preserve $D$. Their restrictions to $D$ are nilpotent and, respectively, symmetric and skew-symmetric with respect to its positive-definite inner product. Hence both restrictions vanish. As both operators also vanish on $D^\perp$, we obtain:
   \[
   \Ll_u = \Rr_u = 0.
   \]

   \item Let \(u \in N(\mathcal{A})\). Then, we have \(\Ll_u = - \Rr_u\), and from \eqref{Mulgauch} we get \(\Ll_u^2 = 0\). Since \(\Ll_u\) is symmetric with respect to the Euclidean metric \(\prs\), it follows that:
   \[
   \Ll_u = \Rr_u = 0.
   \]
\end{enumerate}
\end{proof}

\begin{Proposition}
Let $(\mathcal{A}, \bullet, \prs)$ be a Euclidean left-alternative algebra. Then:
\[
N(\mathcal{A}) = N_{\ell r}(\mathcal{A}) = N_r(\mathcal{A}) = (\mathcal{A} \bullet \mathcal{A})^\perp = (\mathcal{A} \cdot \mathcal{A})^\perp.
\]
\end{Proposition}

\begin{proof}
By \eqref{eq1}, $N_r(\A)=(\A\bullet\A)^\perp\subseteq(\A\cdot\A)^\perp$. Proposition~\ref{Euclidean}(2) yields $(\A\cdot\A)^\perp\subseteq N_{\ell r}(\A)\subseteq N_r(\A)$. Moreover, $N_{\ell r}(\A)\subseteq N(\A)$, and Proposition~\ref{Euclidean}(3) gives the reverse inclusion. Thus all the stated subspaces coincide.
\end{proof}

\begin{Proposition}\label{Trivial algebra}
Let \((\mathcal{A}, \bullet, \prs)\) be a pseudo-Euclidean left-alternative algebra such that the Jordan algebra \(\mathcal{A}^+\) is nilpotent. If \(\mathcal{A} \cdot \mathcal{A}\) is Euclidean, then $\mathcal{A} \bullet \mathcal{A} = \{0\}.$
\end{Proposition}

\begin{proof}
Set $D=\A\cdot\A$. Since the restriction of $\prs$ to $D$ is positive definite, we have the orthogonal decomposition $\A=D\oplus D^\perp$. By Proposition~\ref{Euclidean}(2), $\Ll_u=\Rr_u=0$ for all $u\in D^\perp$. Consequently, $D^\perp\subseteq N(\A)$ and
\[
D=\A\cdot\A=D\cdot D.
\]
Since $D$ is a subalgebra of the nilpotent Jordan algebra $\A^+$, it is nilpotent. The equality $D=D\cdot D$ therefore implies $D=\{0\}$. Finally, by \eqref{kozule}, we obtain $\A\bullet\A=\{0\}$.
\end{proof}

\begin{Corollary}\label{Eutrivial}
Let \((\mathcal{A}, \bullet, \prs)\) be a Euclidean left-alternative algebra such that the Jordan algebra \(\mathcal{A}^+\) is nilpotent. Then, $\mathcal{A} \bullet \mathcal{A} = \{0\}.$
\end{Corollary}

\begin{proof}
Since \(\prs\) is Euclidean on \(\mathcal{A}\), it is also Euclidean on \(\mathcal{A} \cdot \mathcal{A}\). Hence, by Proposition~\ref{Trivial algebra}, $\mathcal{A} \bullet \mathcal{A} = \{0\}.$
\end{proof}

\begin{Proposition}\label{prop de kerL}
Let \((\mathcal{A}, \bullet, \prs)\) be a Euclidean left-alternative algebra, and let \((\mathcal{A}, \cdot)\) denote its associated Jordan algebra. Then:
\begin{enumerate}
    \item The left annihilator \(N_\ell(\mathcal{A})\) is a two-sided ideal of \((\mathcal{A}, \bullet)\).
    
    \item The orthogonal complement $N_\ell(\A)^\perp$, endowed with the restrictions of the product $\bullet$ and of the bilinear form $\prs$, is a Euclidean left-alternative subalgebra. Moreover, its annihilator is trivial, that is:
    \[
    N_{\ell r}(N_\ell(\mathcal{A})^\perp) = \{0\}.
    \]
    
    \item The restriction of the bilinear form \( B_L(u, v) = \operatorname{Tr}(\Ll_u \circ \Ll_v) \) to \(N_\ell(\mathcal{A})^\perp\) is positive-definite and satisfies:
    \[
    B_L(u \cdot v, w) = B_L(u, v \cdot w), \quad \text{for all }\, u, v, w \in N_\ell(\mathcal{A})^\perp.
    \]

    \item For all \( u \in N_\ell(\mathcal{A})^\perp \), let \( \rho(u) \) denote the restriction of \(\operatorname{sd}_u\) to \(N_\ell(\mathcal{A})\). Then, \( \rho(u) \) is symmetric with respect to $\prs$, and for all \( u, v \in N_\ell(\mathcal{A})^\perp \), we have:
    \[
    \rho(u \cdot v) = \rho(u) \circ \rho(v) + \rho(v) \circ \rho(u).
    \]
\end{enumerate}
\end{Proposition}
   
\begin{proof}
\begin{enumerate}
    \item Let \( u \in N_\ell(\mathcal{A}) \) and \( v \in \mathcal{A} \). Then \( u \bullet v = 0 \). Using identity~\eqref{Mulgauch}, we have:
    \[
    \Ll_{v \bullet u} = \Ll_{u \cdot v} = \Ll_u \circ \Ll_v + \Ll_v \circ \Ll_u = 0.
    \]
    Hence, \( v \bullet u \in N_\ell(\mathcal{A}) \), and thus \( N_\ell(\mathcal{A}) \) is a two-sided ideal.

    \item Let \( u, v \in N_\ell(\mathcal{A})^\perp \) and \( w \in N_\ell(\mathcal{A}) \). Then:
    \[
    \langle u \bullet v, w \rangle = \langle v, u \bullet w \rangle = 0,
    \]
    because \( N_\ell(\mathcal{A}) \) is a two-sided ideal. Thus, \( u \bullet v \in N_\ell(\mathcal{A})^\perp \), so \( N_\ell(\mathcal{A})^\perp \) is a subalgebra of \( (\mathcal{A}, \bullet) \). Since \( (N_\ell(\mathcal{A})^\perp, \bullet, \prs) \) is also a Euclidean left-alternative algebra, Proposition~\ref{Euclidean} implies that:
    \[
    N_{\ell r}(N_\ell(\mathcal{A})^\perp) = \{ u \in N_\ell(\mathcal{A})^\perp \mid u \bullet v = v \bullet u = 0, \ \forall v \in N_\ell(\mathcal{A})^\perp \}.
    \]
    In particular, for any \( u \in N_{\ell r}(N_\ell(\mathcal{A})^\perp) \), we have \( u \bullet u = 0 \), so identity~\eqref{Mulgauch} yields \( \Ll_u^2 = 0 \). Since \( \Ll_u \) is symmetric with respect to \(\prs\), this implies \( \Ll_u = 0 \), hence \( u \in N_\ell(\mathcal{A}) \), and thus \( u \in N_\ell(\mathcal{A})^\perp \cap N_\ell(\mathcal{A}) = \{0\} \).

    \item The bilinear form $B_L$ on $\mathcal{A}$ defined by \( B_L(u,v) = \operatorname{Tr}(\Ll_u \circ \Ll_v) \) is symmetric. The fact that $\prs$ is Euclidean and the operator \( \Ll_u \) is symmetric with respect to \(\prs\) for all $u \in \mathcal{A}$ implies that if \( B_L(u,u) = 0 \), then \( \Ll_u = 0 \), which implies \( u \in N_\ell(\mathcal{A}) \). So, the restriction of \( B_L \) to \(N_\ell(\mathcal{A})^\perp\) is positive-definite. Furthermore, using identity~\eqref{Mulgauch}, for all \( u, v, w \in \mathcal{A} \), we obtain:
    \begin{align*}
        B_L(u \cdot v, w) &= \operatorname{Tr}(\Ll_{u \cdot v} \circ \Ll_w) \\
        &= \operatorname{Tr}(\Ll_u \circ \Ll_v \circ \Ll_w) + \operatorname{Tr}(\Ll_v \circ \Ll_u \circ \Ll_w) \\
        &= \operatorname{Tr}(\Ll_u \circ (\Ll_v \circ \Ll_w + \Ll_w \circ \Ll_v)) \\
        &= \operatorname{Tr}(\Ll_u \circ \Ll_{v \cdot w}) = B_L(u, v \cdot w).
    \end{align*}

    \item Let \( u \in N_\ell(\mathcal{A})^\perp \), and define \( \rho(u) := \operatorname{sd}_u\vert_{N_\ell(\mathcal{A})} \). For any \( v \in N_\ell(\mathcal{A})\), we have:
    \[
    \rho(u)(v) = u \cdot v = u \bullet v + v \bullet u = u \bullet v=\Ll_u(v).
    \]
    Since \( \Ll_u \) is symmetric with respect to \(\prs\), so is \(\rho(u)\). Finally, identity~\eqref{Mulgauch} gives:
    \[
    \rho(u \cdot v) = \rho(u) \circ \rho(v) + \rho(v) \circ \rho(u),
    \]
    as desired.
\end{enumerate}
\end{proof}

\begin{Proposition}\label{Semi-simple Ker}
Let $(\mathcal{A}, \bullet, \prs)$ be a Euclidean left-alternative algebra, and let $\mathcal{A}^+$ denote its associated Jordan algebra. Then, the left annihilator $N_\ell( \mathcal{A} )$ coincides with the radical $\mathcal{N}$ of $\mathcal{A}^+$, and its orthogonal complement $N_\ell(\mathcal{A})^\perp$ is a semisimple subalgebra of $\mathcal{A}^+$.
\end{Proposition}

\begin{proof}
Since $\mathcal{A}^+$ is a Jordan algebra, Theorem~\ref{theo decom} ensures the direct sum decomposition $\mathcal{A} = \mathcal{S} \oplus \mathcal{N}$, where $\mathcal{N}$ is the radical of $\mathcal{A}^+$, and $\mathcal{S}$ is a maximal semisimple subalgebra of $\mathcal{A}^+$. Let $u \in \mathcal{N}$. By \eqref{Re7}, we have:
\[
\operatorname{sd}_u^n(u) = 2^n \Ll_u^n(u), \quad \forall n \in \mathbb{N}^*.
\]
Applying Lemma~\ref{Lemme3}, it follows that:
\[
\Ll_{\operatorname{sd}_u^n(u)} = 2^n\Ll_{\Ll_u^n(u)} = 2^n\Ll_u^{n+1}.
\]
As $\mathcal{N}$ is a nilpotent ideal, there exists an integer $n_0$ such that $\operatorname{sd}_u^{n_0}(u) = 0$. Therefore, $\Ll_u^{n_0 + 1} = 0$, which shows that $\Ll_u$ is nilpotent. Since $\Ll_u$ is symmetric and the bilinear form $\prs$ is Euclidean, it follows that $\Ll_u = 0$, hence $u \in N_\ell(\mathcal{A})$. Thus, we obtain $\mathcal{N} \subseteq N_\ell(\mathcal{A})$.

On the other hand, according to Proposition~\ref{prop de kerL}, $N_\ell(\mathcal{A})$ is a two-sided ideal of $(\mathcal{A}, \bullet)$, and hence it is also an ideal of $(\mathcal{A}, \cdot)$. Since $N_\ell(\A)\cdot N_\ell(\A)=\{0\}$, the ideal $N_\ell(\mathcal{A})$ is nilpotent. As $\mathcal{N}$ is the unique maximal nilpotent ideal, it follows that $N_\ell(\mathcal{A}) \subseteq \mathcal{N}$, which yields the equality $N_\ell(\mathcal{A}) = \mathcal{N}$.

Finally, since $N_\ell(\mathcal{A})^\perp$ is isomorphic as a Jordan algebra to the quotient algebra $\mathcal{A} / \mathcal{N}$, it follows that $N_\ell(\mathcal{A})^\perp$ is a semisimple subalgebra of $\mathcal{A}^+$.
\end{proof}

\begin{Definition}
Let $(\mathcal{A}, \bullet)$ be a non-associative algebra. A bilinear form $\prs$ on $\mathcal{A}$ is called \textit{associative} if:
\(
\langle u \bullet v, w \rangle = \langle u, v \bullet w \rangle,
\) 
for all $u, v, w \in \mathcal{A}$. If, in addition, $\prs$ is non-degenerate and symmetric, the triple $(\mathcal{A}, \bullet, \prs)$ is called a \textit{quadratic algebra}.
\end{Definition}

\begin{UnnumberedRemark}
In \cite{Benayadi2} (resp. \cite{Boulmane}), quadratic Jordan algebras (resp. alternative algebras) were called pseudo-Euclidean Jordan algebras (resp. pseudo-Euclidean alternative algebras).
\end{UnnumberedRemark}

The following result, from \cite{Boulmane}, shows that in a quadratic algebra, left and right alternativity are equivalent, and thus the algebra is alternative:

\begin{Proposition}[\cite{Boulmane}]\label{Quadratic}
Let $(\mathcal{A}, \bullet, \prs)$ be a quadratic non-associative algebra. Then, the following assertions are equivalent:
\begin{enumerate}
    \item $(\mathcal{A}, \bullet)$ is a left-alternative algebra.
    \item $(\mathcal{A}, \bullet)$ is a right-alternative algebra.
    \item $(\mathcal{A}, \bullet)$ is an alternative algebra.
\end{enumerate}
\end{Proposition}

\begin{Proposition}\label{element neutre de Levi}
Let $(\mathcal{A}, \bullet, \prs)$ be a Euclidean left-alternative algebra such that its associated Jordan algebra $\mathcal{A}^+$ is semisimple. Then, $e := 2E$ is the unit element of the left-alternative algebra $(\mathcal{A}, \bullet)$, where $E$ is the unit element of the Jordan algebra $\mathcal{A}^+$.
\end{Proposition}

\begin{proof}
Since $\mathcal{A}^+$ is a finite-dimensional semisimple Jordan algebra, it follows from \cite[Theorem 4.7]{Rich} that $\mathcal{A}$ admits a unit element $E$ with respect to the Jordan product $\cdot$ (i.e., $E \cdot u = u$ for all $u \in \mathcal{A}$). In particular:
\[
E \cdot E = E = 2E \bullet E.
\]

On the other hand, identity~\eqref{Mulgauch} implies that for any $u \in \mathcal{A}$:
\[
\Ll_{u \cdot E} = \Ll_u \circ \Ll_E + \Ll_E \circ \Ll_u.
\]
Hence:
\[
\Ll_u = \Ll_u \circ \Ll_E + \Ll_E \circ \Ll_u.
\]
Consequently, we compute:
\begin{align*}
\Ll_E \circ \Ll_u &= \Ll_E \circ \Ll_u \circ \Ll_E + \Ll_E \circ \Ll_E \circ \Ll_u \\
&= \Ll_E \circ \Ll_u \circ \Ll_E + \Ll_{E \bullet E} \circ \Ll_u \\
&= \Ll_E \circ \Ll_u \circ \Ll_E + \tfrac{1}{2} \Ll_E \circ \Ll_u.
\end{align*}
We deduce:
\[
\Ll_E \circ \Ll_u = 2 \Ll_E \circ \Ll_u \circ \Ll_E.
\]
Similarly, one obtains:
\[
\Ll_u \circ \Ll_E = 2 \Ll_E \circ \Ll_u \circ \Ll_E.
\]
Then, $\Ll_E$ commutes with $\Ll_u$. We conclude that:
\[
\Ll_u = 2 \Ll_E \circ \Ll_u = 2 \Ll_u \circ \Ll_E.
\]

Now, for any $u, v \in \mathcal{A}$, we have:
\[
u \bullet v = (u \bullet v) \cdot E = (u \bullet v) \bullet E + E \bullet (u \bullet v),
\]
which implies:
\[
u \bullet v = 2 (u \bullet v) \bullet E = 2 E \bullet (u \bullet v).
\]
Setting $e := 2E$, we obtain:
\[
u \bullet v = e \bullet (u \bullet v) = (u \bullet v) \bullet e.
\]

Finally, since $\mathcal{A}^+$ is semisimple, we have $\mathcal{A} = \mathcal{A} \cdot \mathcal{A}$, and as $\mathcal{A} \cdot \mathcal{A} \subseteq \mathcal{A} \bullet \mathcal{A}$, it follows that $\mathcal{A} = \mathcal{A} \bullet \mathcal{A}$. Thus, for all $u \in \mathcal{A}$:
\[
u = e \bullet u = u \bullet e.
\]
This completes the proof.
\end{proof}

\begin{Proposition}
Let $(\mathcal{A}, \bullet, \prs)$ be a Euclidean left-alternative algebra such that its associated Jordan algebra $\mathcal{A}^+$ is simple. Then, the bilinear forms $B_L$ and $B_R$ on $\mathcal{A}$, defined by
\[
B_L(u,v) = \operatorname{Tr}(\Ll_u \circ \Ll_v), \quad B_R(u,v) = \operatorname{Tr}(\Rr_u \circ \Rr_v) \quad \forall\, u, v \in \mathcal{A},
\]
coincide; that is, $B_L = B_R$.
\end{Proposition}

\begin{proof}
Since $\mathcal{A}^+$ is simple, its radical $\mathcal{N}$ is trivial. Proposition~\ref{Semi-simple Ker} then ensures that $\mathcal{A} = N_\ell(\mathcal{A})^\perp$. By Proposition~\ref{prop de kerL}, the bilinear form $B_L(u,v) = \operatorname{Tr}(\Ll_u \circ \Ll_v)$ is symmetric, associative, and positive-definite on $\mathcal{A}^+$. 

Let $E$ be the unit element of $\A^+$. If $\sum_{j=1}^m u_j\cdot u_j=0$, associativity of $B_L$ yields
\[
0=B_L\left(E,\sum_{j=1}^m u_j\cdot u_j\right)=\sum_{j=1}^m B_L(u_j,u_j).
\]
By positive definiteness, $u_j=0$ for every $j$. Thus $\A^+$ is formally real, and the results on simple Euclidean Jordan algebras in \cite{Faraut} apply.

According to \cite[Proposition~III.1.5]{Faraut}, the bilinear form $\operatorname{Al}(u,v) := \operatorname{Tr}(\operatorname{sd}_{u\cdot v})$ is also symmetric and positive-definite on $\mathcal{A}^+$. Since $\mathcal{A}^+$ is simple, \cite[Proposition~III.4.1]{Faraut} implies that there exists a constant $\lambda > 0$ such that $\operatorname{Al} = \lambda B_L$.

Since $(\mathcal{A}, \bullet)$ is a left-alternative algebra, we can apply the identities \eqref{Mulgauch} and \eqref{Muladroite} to obtain
\[
\operatorname{Al}(u,v) = \operatorname{Tr}(\Ll_{u\cdot v}) + \operatorname{Tr}(\Rr_{u\cdot v}) \overset{\eqref{Mulgauch},\eqref{Muladroite}}{=} 2\operatorname{Tr}(\Ll_u \circ \Ll_v) + 2\operatorname{Tr}(\Rr_u \circ \Rr_v),
\]
which, using the relation $\operatorname{Al} = \lambda B_L$, yields
\[
2\operatorname{Tr}(\Ll_u \circ \Ll_v) + 2\operatorname{Tr}(\Rr_u \circ \Rr_v) = \lambda\operatorname{Tr}(\Ll_u \circ \Ll_v).
\]

By Proposition~\ref{element neutre de Levi}, the algebra $(\mathcal{A}, \bullet)$ admits a unit element $e$, which implies $\Ll_e = \Rr_e = \mathrm{Id}_{\mathcal{A}}$. Thus, 
\[
\operatorname{Tr}(\Ll_e \circ \Ll_e) = \operatorname{Tr}(\Rr_e \circ \Rr_e) = \dim \mathcal{A} = n.
\]
Evaluating the previous relation at $u = v = e$ yields $2n + 2n = \lambda n$, which implies $\lambda = 4$. Consequently, we obtain
\[
\operatorname{Tr}(\Ll_u \circ \Ll_v) = \operatorname{Tr}(\Rr_u \circ \Rr_v)
\]
for all $u, v \in \mathcal{A}$, which completes the proof.
\end{proof}

We now establish that any Euclidean left-alternative algebra whose associated Jordan algebra is simple is an alternative algebra. 

\begin{Proposition}\label{alternative}
Let $(\mathcal{A}, \bullet, \prs)$ be a Euclidean left-alternative algebra such that its associated Jordan algebra $\mathcal{A}^+$ is simple. Then:
\begin{itemize}
\item[(i)] $(\mathcal{A}, \bullet, B_L)$ is a positive-definite quadratic left-alternative algebra;
\item[(ii)] $(\mathcal{A}^+, B_L)$ is a positive-definite quadratic Jordan algebra;
\item[(iii)] $(\mathcal{A}, \bullet)$ is an alternative algebra.
\end{itemize}
\end{Proposition}

\begin{proof}
Since $\mathcal{A}^+$ is simple, its radical $\mathcal{N}$ is trivial. By Proposition~\ref{Semi-simple Ker}, we have $\mathcal{A} = N_\ell(\mathcal{A})^\perp$. According to Proposition~\ref{prop de kerL}, the bilinear form $B_L(u,v) := \operatorname{Tr}(\Ll_u \circ \Ll_v)$ is symmetric, associative, and positive-definite on $\mathcal{A}^+$. For any $u,v,w \in \mathcal{A}$, we have:
\[
B_L(u \cdot v, w) = B_L(u, v \cdot w) = \operatorname{Tr}(\Rr_{u \cdot v} \circ \Rr_w) = \operatorname{Tr}(\Rr_u \circ \Rr_{v \cdot w}).
\]
From identity~\eqref{Muladroite}, we get:
\begin{equation}
\operatorname{Tr}([\Rr_u, \Ll_v] \circ \Rr_w) + \operatorname{Tr}([\Rr_v, \Ll_u] \circ \Rr_w) = \operatorname{Tr}([\Rr_v, \Ll_w] \circ \Rr_u) + \operatorname{Tr}([\Rr_w, \Ll_v] \circ \Rr_u).\label{Re8}
\end{equation}
On the other hand:
\[
\begin{aligned}
\operatorname{Tr}([\Rr_w,\Ll_v]\circ\Rr_u)
&=\operatorname{Tr}(\Rr_w\circ\Ll_v\circ\Rr_u-\Ll_v\circ\Rr_w\circ\Rr_u)\\
&=\operatorname{Tr}(\Ll_v\circ\Rr_u\circ\Rr_w)-\operatorname{Tr}(\Rr_u\circ\Ll_v\circ\Rr_w)\\
&=-\operatorname{Tr}([\Rr_u,\Ll_v]\circ\Rr_w).
\end{aligned}
\]
Inserting this into \eqref{Re8}, we obtain:
\begin{equation}
2 \operatorname{Tr}([\Rr_u, \Ll_v] \circ \Rr_w) = \operatorname{Tr}([\Rr_w, \Ll_u] \circ \Rr_v) + \operatorname{Tr}([\Rr_v, \Ll_w] \circ \Rr_u). \label{Re9}
\end{equation}
On the other hand, we have:
\begin{equation}
\begin{aligned}
2\operatorname{Tr}([\Rr_u, \Ll_v]\circ \Rr_w) &\overset{\eqref{Re9}}{=} \operatorname{Tr}([\Rr_w, \Ll_u]\circ\Rr_v)+\operatorname{Tr}([\Rr_v, \Ll_w]\circ\Rr_u)\\ 
&\overset{\eqref{Re9}}{=} \frac{1}{2}\operatorname{Tr}([\Rr_v, \Ll_w]\circ\Rr_u)+\frac{1}{2}\operatorname{Tr}([\Rr_u, \Ll_v]\circ\Rr_w) +\operatorname{Tr}([\Rr_v, \Ll_w]\circ\Rr_u)\\ 
&= \frac{1}{2}\operatorname{Tr}([\Rr_u, \Ll_v]\circ\Rr_w)+\frac{3}{2}\operatorname{Tr}([\Rr_v, \Ll_w]\circ\Rr_u).
\end{aligned}
\end{equation}
Thus, we derive the cyclic identity:
\begin{equation}
\operatorname{Tr}([\Rr_u, \Ll_v] \circ \Rr_w) = \operatorname{Tr}([\Rr_v, \Ll_w] \circ \Rr_u) = \operatorname{Tr}([\Rr_w, \Ll_u] \circ \Rr_v).
\label{cyclic} \end{equation}
Now, for any $u,v,w \in \mathcal{A}$, we have:
\[
\begin{aligned}
B_L(u \bullet v, w) &= \operatorname{Tr}(\Rr_{u \bullet v} \circ \Rr_w) \overset{\eqref{Muladroite}}{=} \operatorname{Tr}(\Rr_v \circ \Rr_u \circ \Rr_w + [\Rr_v, \Ll_u] \circ \Rr_w) \\
&\overset{\eqref{cyclic}}{=} \operatorname{Tr}(\Rr_w \circ \Rr_v \circ \Rr_u) + \operatorname{Tr}([\Rr_w, \Ll_v] \circ \Rr_u) \\ 
&= \operatorname{Tr}((\Rr_w \circ \Rr_v + [\Rr_w, \Ll_v]) \circ \Rr_u) \\
&\overset{\eqref{Muladroite}}{=} \operatorname{Tr}(\Rr_u \circ \Rr_{v \bullet w}) = B_L(u, v \bullet w).
\end{aligned}
\]
Thus, the bilinear form $B_L$ is associative with respect to the product $\bullet$, and hence $(\mathcal{A}, \bullet, B_L)$ is a quadratic left-alternative algebra. By Proposition~\ref{Quadratic}, it follows that $(\mathcal{A}, \bullet)$ is an alternative algebra.
\end{proof}

\section{Characterization of Euclidean Left-Alternative Algebras}

In this section, we first show that if a Euclidean left-alternative algebra has a simple associated Jordan algebra $\mathcal{A}^+ = (\mathcal{A}, \cdot)$, then it must be one-dimensional. We then provide a complete characterization of Euclidean left-alternative algebras. In particular, we prove that every Euclidean left-alternative algebra is associative.

Now, we recall some important tools that will be needed in the results of this section. 

A Jordan algebra $(\mathcal{A}, \cdot)$ over $\mathbb{R}$ with an identity element $E$ is called a \textit{formally real Jordan algebra} if:
\[
u^2 + v^2 = 0 \implies u = 0 \text{ and } v = 0.
\]

Now, let $(\mathcal{A}, \cdot)$ be a formally real Jordan algebra with unit element $E$, and let $\mathbb{R}[Y]$ denote the subalgebra of $\mathcal{A}$ generated by an element $Y \in \mathcal{A}$ and the identity element. Since $\mathbb{R}[Y]$ is finite-dimensional, the number $\operatorname{dim} \mathbb{R}[Y]$ is clearly bounded. We define the rank $r$ of $\mathcal{A}$ as follows:
\[
r = \max \{\operatorname{dim} \mathbb{R}[Y] \mid Y \in \mathcal{A}\}.
\]
An element $X \in \mathcal{A}$ is said to be \textit{regular} if $\operatorname{dim} \mathbb{R}[X]$ is maximal among $\operatorname{dim} \mathbb{R}[Y]$ for all $Y \in \mathcal{A}$. The reduced trace of a regular element $X \in \mathcal{A}$ is defined by:
\begin{equation}\label{Re3}
\operatorname{tr}(X) = \left. \operatorname{Tr} \operatorname{sd}_X \right|_{\mathbb{R}[X]}.
\end{equation}
Moreover, since the set of regular elements is open and dense in $\mathcal{A}$, the function \eqref{Re3} can be uniquely extended to a linear map $\operatorname{tr}: \mathcal{A} \rightarrow \mathbb{R}$, and in particular, we have $\operatorname{tr}(E) = r$ (see \cite{Faraut}).

A non-zero element $C$ of $\mathcal{A}$ is called a \textit{primitive idempotent} if $C^2 = C$ and $C$ cannot be written as the sum of two nonzero orthogonal idempotents. It is known that there exists a \textit{Jordan frame}, i.e., a set $\{H_1, H_2, \dots, H_r\}$ of primitive idempotents of $(\mathcal{A}, \cdot)$ such that $H_1 + H_2 + \dots + H_r$ is the identity element $E$, and $H_i \cdot H_j = 0$ for all $i \neq j$. In a simple formally real Jordan algebra, Jordan frames are unique up to automorphisms (and permutation of their elements).

\begin{Theorem}[Spectral theorem, second version \cite{Faraut}, p. 44]\label{Spectral}
Let $(\mathcal{A}, \cdot)$ be a finite-dimensional formally real Jordan algebra with unit element $E$ and rank $r$. Then, for any $u \in \mathcal{A}$, there exist a Jordan frame $\{H_1, \dots, H_r\}$ and real numbers $\lambda_1, \dots, \lambda_r$ such that:
\[
u = \sum_{j=1}^r \lambda_j H_j,
\]
where the numbers $\lambda_1, \dots, \lambda_r \in \mathbb{R}$ are uniquely determined by $u$ up to permutation. In this case, we have $\operatorname{tr}(u) = \lambda_1 + \dots + \lambda_r$. In particular, $\operatorname{tr}(H_i) = 1$ for all $1 \leq i \leq r$, and $\operatorname{tr}(E) = r$.
\end{Theorem}

\begin{Remark}
On a simple Euclidean Jordan algebra, any two associative inner products are proportional. Consequently, orthogonality is independent of the choice of the associative inner product.
\end{Remark}

\begin{Theorem}[\cite{Faraut}, p. 68] \label{Faraut2}
Let $(\mathcal{A}, \cdot)$ be a simple formally real Jordan algebra of rank $r$ and let $\{H_1, H_2, \dots, H_r\}$ be a Jordan frame of $(\mathcal{A}, \cdot)$. Then, $\mathcal{A}$ decomposes into the following orthogonal direct sum (with respect to an associative inner product):
\[
\mathcal{A} = \bigoplus_{1 \leq i \leq j \leq r} \mathcal{A}_{i j} = \big( \bigoplus_{i=1}^r \mathbb{R} H_i \big) \oplus \big( \bigoplus_{1 \leq i < j \leq r} \mathcal{A}_{i j} \big).
\]
Moreover, we have:
\[
\begin{aligned}
&\mathcal{A}_{i j}\cdot\mathcal{A}_{i j} \subseteq \mathcal{A}_{i i} + \mathcal{A}_{j j}, \quad i\neq j, \\
&\mathcal{A}_{i j}\cdot\mathcal{A}_{j k} \subseteq \mathcal{A}_{i k}, \quad \text{if }i,j,k\text{ are pairwise distinct}, \\
&\mathcal{A}_{i j} \mathcal{A}_{k l} = \{0\}, \quad \text{if } \{i, j\} \cap \{k, l\} = \varnothing,
\end{aligned}
\]
where $\mathcal{A}_{ii}=\mathbb{R}H_i$, and for $i<j$,
\[
\mathcal{A}_{ij}=\{u\in\A\mid H_i\cdot u=H_j\cdot u=\tfrac12u\},
\qquad \mathcal{A}_{ji}=\mathcal{A}_{ij}.
\]
In particular, $H_k\cdot u=0$ for $u\in\mathcal{A}_{ij}$ and $k\notin\{i,j\}$.
\end{Theorem}

In \cite{Faraut}, it is shown that if $(\mathcal{A}, \cdot, B)$ is a quadratic Jordan algebra with a positive-definite bilinear form $B$ and possessing a unit element $E$, then $(\mathcal{A}, \cdot)$ is formally real. Indeed, for any $u, v \in \mathcal{A}$ such that $u^2 + v^2 = 0$, we have:
\[
B(u, u) + B(v, v) = B(E \cdot u, u) + B(E \cdot v, v) = B(E, u^2) + B(E, v^2) = B(E, u^2 + v^2) = 0.
\]
By the positive definiteness of $B$, it follows that $u = v = 0$. Hence, $\mathcal{A}$ is formally real.

We now establish that this result applies in the context of certain Euclidean left-alternative algebras:

\begin{Proposition}\label{formally}
Let $(\mathcal{A}, \bullet, \prs)$ be a Euclidean left-alternative algebra such that its associated Jordan algebra $\mathcal{A}^+ = (\mathcal{A}, \cdot)$ is simple. Then, $\mathcal{A}^+$ is a formally real Jordan algebra.
\end{Proposition}

\begin{proof}
By Propositions~\ref{prop de kerL} and~\ref{Semi-simple Ker}, $(\mathcal{A}^+, B_L)$ is a quadratic Jordan algebra with a positive-definite associative bilinear form $B_L$. Applying the result from \cite{Faraut}, it follows that $\mathcal{A}^+$ is formally real.
\end{proof}

\begin{Lemma}\label{Lemme6}
Let $(\mathcal{A}, \bullet, \prs)$ be a Euclidean left-alternative algebra such that its associated Jordan algebra $\mathcal{A}^+ = (\mathcal{A}, \cdot)$ is simple. Let $\{H_1, H_2, \dots, H_r\}$ be a Jordan frame of $(\mathcal{A}, \cdot)$. Then, $\mathcal{A}$ splits orthogonally as:
\begin{equation}
\mathcal{A} = \bigoplus_{1 \leq i \leq j \leq r} \mathcal{A}_{i j} = \big( \bigoplus_{i=1}^r \mathbb{R} H_i \big) \oplus \big( \bigoplus_{1 \leq i < j \leq r} \mathcal{A}_{i j} \big).
\end{equation}
\end{Lemma}

\begin{proof}
By Proposition~\ref{formally}, the Jordan algebra $\mathcal{A}^+$ is formally real. Hence, we can apply Theorem~\ref{Faraut2}. In particular, for any $i \in \{1, \dots, r\}$, we have:
\[
H_i \cdot H_i = H_i = 2 H_i \bullet H_i.
\]
If $H_i \cdot u = 0$ for $u \in \mathcal{A}$, then it follows that $H_i \bullet u = - u \bullet H_i$. We are going to show that $H_i \bullet u = u \bullet H_i = 0$. Since $(\mathcal{A}, \bullet, \prs)$ is a Euclidean left-alternative algebra and its associated Jordan algebra $\mathcal{A}^+$ is simple, Proposition~\ref{alternative} implies that $(\mathcal{A}, \bullet)$ is alternative. So, we have:
\begin{equation}
(H_i \bullet u) \bullet H_i = H_i \bullet (u \bullet H_i), \quad
H_i \bullet (H_i \bullet u) = (H_i \bullet H_i) \bullet u, \quad
u \bullet (H_i \bullet H_i) = (u \bullet H_i) \bullet H_i.
\label{u0}\end{equation}
Then:
\[
\begin{aligned}
H_i \bullet u &= 2 (H_i \bullet H_i) \bullet u = 2 H_i \bullet (H_i \bullet u) = -2 H_i \bullet (u \bullet H_i) \\
&= -2 (H_i \bullet u) \bullet H_i = 2 (u \bullet H_i) \bullet H_i = 2 u \bullet (H_i \bullet H_i) = u \bullet H_i.
\end{aligned}
\]
Consequently, $H_i\bullet u=u\bullet H_i$. Together with $H_i\bullet u=-u\bullet H_i$, this gives:
\begin{equation}\label{L_0}
H_i \bullet u = u \bullet H_i = 0.
\end{equation}
Consequently:
\[
\langle H_i, u \rangle = 2 \langle H_i \bullet H_i, u \rangle = 2 \langle H_i, H_i \bullet u \rangle = 0.
\]
Let $i\neq j$ and let $k<l$ with $i\notin\{k,l\}$. For any $s_{kl}\in\mathcal{A}_{kl}$, since $H_i \cdot H_j = 0$, $H_i \cdot s_{kl} = 0$, and $H_i \cdot H_i = 2 H_i \bullet H_i$, we also have:
\[
\langle H_i, H_j \rangle = \langle H_i, s_{kl} \rangle = 0.
\]

Recall that $E = \sum_{i=1}^r H_i$ is the unity of $\mathcal{A}^+$. By Proposition~\ref{element neutre de Levi}, it follows that:
\[
2 E \bullet u = 2 u \bullet E = u, \quad \text{ for all } u \in \mathcal{A}.
\]
For $s_{ij} \in \mathcal{A}_{ij}$, from the above relations, we deduce:
\begin{equation}\label{L1}
s_{ij} = 2 s_{ij} \bullet E = 2 E \bullet s_{ij} = 2 (H_i + H_j) \bullet s_{ij} = 2 s_{ij} \bullet (H_i + H_j).
\end{equation}

Next, we are going to prove that $\langle H_i, s_{ij} \rangle = 0$. Using Koszul's formula \eqref{kozule}:
\[
\begin{aligned}
2 \langle s_{ij} \bullet H_i, H_i \rangle &= \langle s_{ij} \cdot H_i, H_i \rangle - \langle H_i \cdot H_i, s_{ij} \rangle + \langle H_i \cdot s_{ij}, H_i \rangle \\
&= \frac{1}{2} \langle s_{ij}, H_i \rangle - \langle H_i, s_{ij} \rangle + \frac{1}{2} \langle s_{ij}, H_i \rangle = 0.\;\;   (\star) \quad 
\end{aligned}
\]
Since $(\mathcal{A}, \bullet)$ is alternative and $H_i \bullet H_j = 0$, we have:
\begin{equation}\label{L2}
H_i \bullet (s_{ij} \bullet H_j) = (H_i \bullet s_{ij}) \bullet H_j - H_i \bullet (H_j \bullet s_{ij}) + (H_i \bullet H_j) \bullet s_{ij} = (H_i \bullet s_{ij}) \bullet H_j - H_i \bullet (H_j \bullet s_{ij}).
\end{equation}

Hence:
\[
\begin{aligned}
\langle H_i, s_{ij} \rangle &\overset{\eqref{L1}}{=} 2 \langle H_i, s_{ij} \bullet (H_i + H_j) \rangle \overset{(\star)}{=} 2 \langle H_i, s_{ij} \bullet H_j \rangle \\
&= 4 \langle H_i \bullet H_i, s_{ij} \bullet H_j \rangle = 4 \langle H_i, H_i \bullet (s_{ij} \bullet H_j) \rangle \\
&\overset{\eqref{L2}}{=} 4 \langle H_i, (H_i \bullet s_{ij}) \bullet H_j - H_i \bullet (H_j \bullet s_{ij}) \rangle \\
&= 4 \langle (H_i \bullet s_{ij}) \bullet H_i, H_j \rangle - 4 \langle H_i \bullet H_i, H_j \bullet s_{ij} \rangle \\
&\overset{\eqref{u0}}{=} 4 \langle H_i \bullet (s_{ij} \bullet H_i), H_j \rangle - 2 \langle H_i, H_j \bullet s_{ij} \rangle \\
&= 4 \langle s_{ij} \bullet H_i, H_i \bullet H_j \rangle - 2 \langle H_j \bullet H_i, s_{ij} \rangle=0.
\end{aligned}
\]

Similarly, for $s_{ij} \in \mathcal{A}_{ij}$ and $s_{jk} \in \mathcal{A}_{jk}$ with $k \notin \{i,j\}$, where $i, j, k \in \{1, \dots, r\}$, we are going to prove that:
\[
\langle s_{ij}, s_{jk} \rangle = 0.
\]
Indeed:
\[
\begin{aligned}
\langle s_{ij}, s_{jk} \rangle &\overset{\eqref{L1}}{=} 2 \langle s_{ij}, (H_j + H_k) \bullet s_{jk} \rangle \\
&= 2 \langle s_{ij}, H_j \bullet s_{jk} \rangle + 2 \langle H_k \bullet s_{ij}, s_{jk} \rangle \\
&\overset{\eqref{L_0}}{=} 2 \langle s_{ij}, H_j \bullet s_{jk} \rangle \\
&= 2 \langle s_{ij}, s_{jk} \cdot H_j \rangle - 2 \langle s_{ij}, s_{jk} \bullet H_j \rangle \\
&= \langle s_{ij}, s_{jk} \rangle - 2 \langle s_{ij}, s_{jk} \bullet H_j \rangle.
\end{aligned}
\]
It follows that $\langle s_{ij}, s_{jk} \bullet H_j \rangle = 0$ and:
\[
\begin{aligned}
\langle s_{ij}, s_{jk} \rangle &\overset{\eqref{L1}}{=} 2 \langle s_{ij}, s_{jk} \bullet (H_j + H_k) \rangle = 2 \langle s_{ij}, s_{jk} \bullet H_k \rangle.
\end{aligned}
\]
Next:
\[
\begin{aligned}
\langle s_{ij}, s_{jk} \rangle &\overset{\eqref{kozule}}{=} \langle s_{jk} \cdot H_k, s_{ij} \rangle - \langle s_{ij} \cdot H_k, s_{jk} \rangle + \langle s_{ij} \cdot s_{jk}, H_k \rangle \\
&= \frac{1}{2} \langle s_{jk}, s_{ij} \rangle + \langle s_{ij} \cdot s_{jk}, H_k \rangle.
\end{aligned}
\]
Since $s_{ij} \cdot s_{jk} \in \mathcal{A}_{ik}$ and $\langle \mathcal{A}_{ik}, H_k \rangle =\{ 0\}$, it follows that:
\[
\langle s_{ij}, s_{jk} \rangle = 0.
\]

Finally, for $s_{ij} \in \mathcal{A}_{ij}$ and $s_{kl} \in \mathcal{A}_{kl}$ with disjoint index sets:
\[
\langle s_{ij}, s_{kl} \rangle \overset{\eqref{L1}}{=} 2 \langle (H_i + H_j) \bullet s_{ij}, s_{kl} \rangle = 2 \langle s_{ij}, (H_i + H_j) \bullet s_{kl} \rangle = 0,
\]
which completes the proof.
\end{proof}

\begin{Lemma}\label{Lemme5}
Let $(\mathcal{A}, \bullet, \prs)$ be a Euclidean left-alternative algebra such that its associated Jordan algebra $\mathcal{A}^+ = (\mathcal{A}, \cdot)$ is simple. Then:
\begin{equation}\label{Re12}
\langle u \cdot (v \cdot v), u \rangle = \langle u \cdot u, v \cdot v \rangle, \quad \text{for all } u,v \in \mathcal{A}.
\end{equation}
\end{Lemma}

\begin{proof}
Since $\mathcal{A}^+$ is a simple Jordan algebra, it admits a unit element $E$. By the spectral theorem (Theorem~\ref{Spectral}), for any element $u \in \mathcal{A}$, there exist a Jordan frame $\{H_1, \dots, H_r\}$ and real scalars $\lambda_1, \dots, \lambda_r \in \mathbb{R}$ such that:
\[
u = \sum_{i=1}^r \lambda_i H_i.
\]
Using the Peirce decomposition of $\mathcal{A}$ with respect to the Jordan frame $\{H_i\}$ and Lemma~\ref{Lemme6}, we have the orthogonal decomposition with respect to $\prs$:
\[
\mathcal{A} = \bigoplus_{1 \leq i \leq j \leq r} \mathcal{A}_{i j} = \big( \bigoplus_{i=1}^r \mathbb{R} H_i \big) \oplus \big( \bigoplus_{1 \leq i < j \leq r} \mathcal{A}_{i j} \big).
\]
For any $v \in \mathcal{A}$, we can write:
\[
v \cdot v = \sum_{i=1}^r \mu_i H_i + \sum_{1 \leq i < j \leq r} s_{ij},
\]
where $\mu_i \in \mathbb{R}$ and $s_{ij} \in \mathcal{A}_{ij}$.

Then, by the Peirce decomposition of $\mathcal{A}^+$:
\[
u \cdot (v \cdot v) = \sum_{i=1}^r \lambda_i \mu_i H_i + \sum_{1 \leq i < j \leq r} \frac{\lambda_i + \lambda_j}{2} s_{ij}.
\]
Hence:
\[
\langle u \cdot (v \cdot v), u \rangle = \sum_{i=1}^r \lambda_i^2 \mu_i \langle H_i, H_i \rangle,
\]
and:
\[
\langle u \cdot u, v \cdot v \rangle = \sum_{i=1}^r \lambda_i^2 \mu_i \langle H_i, H_i \rangle.
\]
Therefore:
\[
\langle u \cdot (v \cdot v), u \rangle = \langle u \cdot u, v \cdot v \rangle,
\]
as desired.
\end{proof}

In what follows, we will prove that the simplicity of the Jordan algebra $\mathcal{A}^+$ forces the Euclidean left-alternative algebra $(\mathcal{A}, \bullet, \prs)$ to be one-dimensional.

\begin{Proposition}\label{Simple dim 1}
Let $(\mathcal{A}, \bullet, \prs)$ be a Euclidean left-alternative algebra such that its associated Jordan algebra $\mathcal{A}^+ = (\mathcal{A}, \cdot)$ is simple. Then, $\dim \mathcal{A} = 1$.
\end{Proposition}

\begin{proof}
For any $u, v \in \mathcal{A}$, the alternative identities yield:
\[
u \bullet (u \bullet v) = (u \bullet u) \bullet v, \quad
v \bullet (u \bullet u) = (v \bullet u) \bullet u, \quad
u \bullet (v \bullet u) = (u \bullet v) \bullet u.
\]
We use these identities to compute:
\[
\begin{aligned}
\langle u \cdot v, u \cdot v \rangle
&= \langle u \bullet v + v \bullet u, u \bullet v + v \bullet u \rangle \\
&= \langle u \bullet v, u \bullet v \rangle + 2 \langle u \bullet v, v \bullet u \rangle + \langle v \bullet u, v \bullet u \rangle \\
&= \langle v, u \bullet (u \bullet v) \rangle + 2 \langle v, u \bullet (v \bullet u) \rangle + \langle u, v \bullet (v \bullet u) \rangle.
\end{aligned}
\]
From Proposition~\ref{alternative}, we know that $(\mathcal{A}, \bullet)$ is an alternative algebra, so:
\[
\begin{aligned}
\langle u \cdot v, u \cdot v \rangle
&= \langle v, (u \bullet u) \bullet v \rangle + 2 \langle v, (u \bullet v) \bullet u \rangle + \langle u, (v \bullet v) \bullet u \rangle \\
&= \langle (u \bullet u) \bullet v, v \rangle + 2 \langle (u \bullet v) \bullet u, v \rangle + \langle (v \bullet v) \bullet u, u \rangle \\
&= \langle (u \bullet u) \bullet v, v \rangle + 2 \langle u \bullet u, v \bullet v \rangle + \langle (v \bullet v) \bullet u, u \rangle \\
&= \langle (u \bullet u) \bullet v, v \rangle + \langle v \bullet (u \bullet u), v \rangle + \langle u, u \bullet (v \bullet v) \rangle + \langle (v \bullet v) \bullet u, u \rangle \\
&= \frac{1}{2} \langle v \cdot (u \cdot u) , v \rangle + \frac{1}{2} \langle u \cdot (v \cdot v), u \rangle.
\end{aligned}
\]
From Lemma~\ref{Lemme5}, we also have:
\begin{equation}\label{Re13}
\langle u \cdot v, u \cdot v \rangle = \langle u \cdot u, v \cdot v \rangle = \langle (u \cdot u) \cdot v, v \rangle.
\end{equation}

Let $\{H_1,\ldots,H_r\}$ be a Jordan frame of $\A^+$. For $1\leq i<j\leq r$ and $v\in\A_{ij}$, we have $H_i\cdot v=\tfrac12v$ and $H_i\cdot H_i=H_i$. Applying \eqref{Re13} with $u=H_i$, we obtain
\[
\begin{aligned}
\tfrac14\langle v,v\rangle
&=\langle H_i\cdot v,H_i\cdot v\rangle\\
&=\langle(H_i\cdot H_i)\cdot v,v\rangle
=\tfrac12\langle v,v\rangle.
\end{aligned}
\]
Since $\prs$ is positive definite, $v=0$. Consequently, $\A_{ij}=\{0\}$ for every $i<j$. Therefore, the Peirce decomposition reduces to
\[
\A=\bigoplus_{i=1}^r\mathbb{R}H_i.
\]
Since $H_i\cdot H_j=0$ for $i\neq j$, each $\mathbb{R}H_i$ is a nonzero ideal of $\A^+$. The simplicity of $\A^+$ implies that $r=1$, and hence $\dim\A=1$.
\end{proof}

Let us recall that if $(\mathcal{A}, \cdot)$ is a non-associative algebra, an element of $\mathcal{A}$ is called an \textit{idempotent} if one has $e \cdot e = e$. A set $\{e_1, \dots, e_p\}$ of elements of $\mathcal{A}$ is called a set of \textit{orthogonal idempotents} of $\mathcal{A}$ if for all $i \in \{1, \dots, p\}$, $e_i$ is an idempotent and for all $i, j \in \{1, \dots, p\}$ such that $i \neq j$, one has $e_i \cdot e_j = 0$.

We are now in a position to state the principal result of this article, which provides a complete structural characterization of Euclidean left-alternative algebras.

\begin{Theorem}\label{Theo Milnor}
Let $(\mathcal{A}, \bullet, \prs)$ be a non-trivial Euclidean left-alternative algebra. Then, there exist pairwise orthogonal idempotents $\{e_1, \dots, e_n\}$ of $(\mathcal{A}, \bullet)$ and an integer $0\leq p\leq n$ such that the following assertions hold:
\begin{enumerate}
\item For each $i, j \in \{1, \dots, n\}$, one has:
\[
\langle e_i,e_i\rangle=\lambda_i>0,\qquad \langle e_i,e_j\rangle=0\quad\text{for }i\neq j.
\]

\item For each $i \in \{1, \dots, p\}$, the operator $\rho_i := \Ll_{e_i}\vert_{N_\ell(\mathcal{A})}$ is diagonalizable and has eigenvalues contained in $\{0, 1\}$.

\item For each $i \in \{1, \dots, p\}$, let $\mathcal{B}_i$ denote the eigenspace of $\rho_i$ associated with the eigenvalue $1$. Then, $\mathcal{B}_i$ is a non-degenerate two-sided ideal of $(\mathcal{A}, \bullet, \prs)$.

\item For all $i, j \in \{1, \dots, p\}$ with $i \neq j$, the ideals $\mathcal{B}_i$ and $\mathcal{B}_j$ are orthogonal, that is:
\[
\langle \mathcal{B}_i, \mathcal{B}_j \rangle = \{0\}.
\]

\item For all $1\leq i\leq p$ and $1\leq j\leq n$, one has $\langle N_{\ell r}(\A),\mathcal{B}_i\rangle=\langle e_j,\mathcal{B}_i\rangle=\{0\}$. Moreover, $\langle e_j,N_{\ell r}(\A)\rangle=\{0\}$ for every $1\leq j\leq n$.

\item The algebra $\mathcal{A}$ admits the orthogonal direct sum decomposition:
\[
\mathcal{A} = N_{\ell r}(\mathcal{A}) \;\oplus\; \bigoplus_{i=1}^p \mathcal{B}_i \;\oplus\; \operatorname{span}\{e_1, \dots, e_n\},
\]
where:
\[
N_\ell(\mathcal{A}) = N_{\ell r}(\mathcal{A}) \;\oplus\; \bigoplus_{i=1}^p \mathcal{B}_i.
\]
\end{enumerate}
\end{Theorem}

\begin{proof}
Since the scalar product $\prs$ is Euclidean, we can decompose the algebra $\mathcal{A}$ orthogonally as:
\[
\mathcal{A} = N_{\ell}(\mathcal{A}) \oplus N_{\ell}(\mathcal{A})^\perp.
\]
According to Proposition~\ref{Semi-simple Ker}, the orthogonal complement $N_{\ell}(\mathcal{A})^\perp$ is a semisimple subalgebra of the Jordan algebra $\mathcal{A}^+$. Hence, there exists a decomposition:
\[
N_{\ell}(\mathcal{A})^\perp = \mathcal{A}_1 \oplus \dots \oplus \mathcal{A}_n,
\]
where each $\mathcal{A}_i$ is a simple ideal of $N_{\ell}(\mathcal{A})^\perp$. We first show that these ideals are pairwise orthogonal. Since $\mathcal{A}_i$ is simple, it admits a unit element $E_i$. Now, take $i \neq j$ and let $u \in \mathcal{A}_j$. Then, $E_i \cdot u = 0$, which implies $E_i \bullet u = - u \bullet E_i$. From the left-alternative identity, we have:
\[
u \bullet E_i = 2 u \bullet (E_i \bullet E_i) = -2 E_i \bullet (u \bullet E_i) = 2 E_i \bullet (E_i \bullet u) = 2 (E_i \bullet E_i) \bullet u = E_i \bullet u.
\]
Hence, we deduce that:
\[
E_i \bullet u = u \bullet E_i = 0.
\]
Let $E = E_1 + \dots + E_n$; this element is the unit of $N_{\ell}(\mathcal{A})^\perp$. Then, according to Proposition~\ref{prop de kerL}, the restriction of $\bullet$ to $N_{\ell}(\mathcal{A})^\perp$ defines a left-alternative algebra. From Proposition~\ref{element neutre de Levi}, we deduce that for all $u \in N_{\ell}(\mathcal{A})^\perp$, we have $2 E \bullet u = 2 u \bullet E = u$.

For $u\in\A_i$, the preceding vanishing identities imply $u=2E\bullet u=2E_i\bullet u$. Now, take $u\in\A_i$ and $v\in\A_j$ with $i\neq j$. Then:
\[
\langle u, v \rangle = 2 \langle E_i \bullet u, v \rangle = 2 \langle u, E_i \bullet v \rangle = 0,
\]
so the ideals $\mathcal{A}_i$ are pairwise orthogonal.

Next, we show that each $\A_i$ is stable under the product $\bullet$. Set $S=N_\ell(\A)^\perp$. Let $u,v\in\A_i$ and $w\in\bigoplus_{j\neq i}\A_j$. Since these are pairwise orthogonal ideals of the Jordan algebra $(S,\cdot)$, we have $u\cdot v\in\A_i$ and $u\cdot w=v\cdot w=0$. By the Koszul formula \eqref{kozule},
\[
2\langle u\bullet v,w\rangle=\langle u\cdot v,w\rangle-\langle v\cdot w,u\rangle+\langle u\cdot w,v\rangle=0.
\]
Since $S$ is a subalgebra of $(\A,\bullet)$, we also have $u\bullet v\in S$. Hence $u\bullet v\in\A_i$, so each $\A_i$ is a Euclidean left-alternative subalgebra whose associated Jordan algebra is simple.

From Proposition~\ref{Simple dim 1}, each simple ideal $\mathcal{A}_i$ is one-dimensional. Set $e_i=2E_i\in\A_i$. Then:
\[
e_i \bullet e_i = e_i \quad \text{and} \quad \langle e_i, e_i \rangle = \lambda_i,
\]
where $\lambda_i>0$. Applying the Koszul formula once again, we obtain $e_i \bullet e_j = 0$ with $i \neq j$.

From the left-alternativity identity \eqref{Mulgauch}, we have:
\[
\Ll_{e_i} = \Ll_{e_i}^2 \quad \text{and} \quad \Ll_{e_i} \circ \Ll_{e_j} = -\Ll_{e_j} \circ \Ll_{e_i}, \quad i \neq j.
\]
Let $\rho_i := \Ll_{e_i}\vert_{N_\ell(\mathcal{A})}$. Then, $\rho_i^2 = \rho_i$ and $\rho_i \circ \rho_j = -\rho_j \circ \rho_i$ for $i \neq j$. Consequently, the spectrum of $\rho_i$ is contained in $\{0, 1\}$.

After reordering the idempotents, we may assume that $\rho_1,\ldots,\rho_p$ are nonzero and $\rho_{p+1},\ldots,\rho_n$ vanish, where $0\leq p\leq n$. For $i\neq j$, multiplying the anticommutativity relation on the left and on the right by $\rho_i$ gives
\[
\rho_i\rho_j=-\rho_i\rho_j\rho_i,
\qquad \rho_i\rho_j\rho_i=-\rho_j\rho_i.
\]
Thus $\rho_i\rho_j=\rho_j\rho_i$ and, by anticommutativity, $\rho_i\rho_j=\rho_j\rho_i=0$. The operators $\rho_i$ are symmetric because the $\Ll_{e_i}$ are symmetric and preserve $N_\ell(\A)$. Hence they are orthogonal projections. For $1\leq i\leq p$, put $\mathcal{B}_i=\operatorname{Im}\rho_i=\ker(\rho_i-I)$. Their images are pairwise orthogonal, since for $b_i\in\mathcal{B}_i$ and $b_j\in\mathcal{B}_j$ with $i\neq j$,
\[
\langle b_i,b_j\rangle=\langle\rho_i b_i,b_j\rangle=\langle b_i,\rho_i b_j\rangle=0.
\]
In addition, $\rho=\sum_{i=1}^p\rho_i$ is an orthogonal projection with image $\bigoplus_{i=1}^p\mathcal{B}_i$ and kernel $\mathcal{B}_0=\bigcap_{i=1}^p\ker\rho_i$ (the empty intersection is $N_\ell(\A)$). Indeed, $\rho_j\rho=\rho_j$ for every $j$, so $\ker\rho=\mathcal{B}_0$. Consequently,
\[
N_\ell(\A)=\mathcal{B}_0\oplus\bigoplus_{i=1}^p\mathcal{B}_i
\]
is an orthogonal direct sum. For $b\in N_\ell(\A)$ and $x=a+\sum_{j=1}^n\alpha_j e_j$, where $a\in N_\ell(\A)$, we have $b\bullet x=0$ and $x\bullet b=\sum_{j=1}^n\alpha_j\rho_j(b)$. Thus $\mathcal{B}_0=N_{\ell r}(\A)$. Moreover, if $b\in\mathcal{B}_i$, then $x\bullet b=\alpha_i b\in\mathcal{B}_i$ and $b\bullet x=0$. Therefore each $\mathcal{B}_i$ is a two-sided ideal, and it is non-degenerate by positive definiteness. Since $e_j\in N_\ell(\A)^\perp$ for every $j$, the required orthogonality and the asserted decomposition of $\A$ follow.
\end{proof}

\begin{Corollary}\label{co1}
Let $(\A,\bullet,\prs)$ be a non-trivial Euclidean left-alternative algebra. 
Then there exists a basis of $\A$ given by:
$$
\{e_1,\ldots,e_n,\;\mu_1,\ldots,\mu_s,\; 
f_{11},\ldots,f_{1t_1},\;
\ldots,\;
f_{p1},\ldots,f_{pt_p}\}
$$
where $n\in\mathbb{N}^*$, $0\leq p\leq n$, $s\in\mathbb{N}$, and $t_1,\ldots,t_p\in\mathbb{N}^*$ (with the convention that any family of vectors is empty if its corresponding upper index is $0$). 
With respect to this basis, the only nonzero products $\bullet$ and the only nonzero values of the bilinear form $\prs$ occur among the following:
\[
e_i \bullet e_i = e_i, \qquad
e_k \bullet f_{k t} = f_{kt}, \qquad
\langle e_i, e_i\rangle = \lambda_i, \qquad
\langle f_{k t}, f_{k t}\rangle = 1, \qquad
\langle \mu_r, \mu_r\rangle = 1,
\]
for all $i\in \{1, \cdots, n\}$, $k\in\{1,\ldots, p\}$, $t\in\{1,\ldots, t_k\}$, and $r\in \{1, \cdots, s\}$, where $\la_i >0$.
\end{Corollary}

\begin{proof}
The result follows directly from Theorem~\ref{Theo Milnor}.
\end{proof}

\begin{Corollary}\label{associative}
Let $(\mathcal{A}, \bullet, \prs)$ be a Euclidean left-alternative algebra. Then, $(\mathcal{A}, \bullet)$ is an associative algebra.
\end{Corollary}

\begin{proof}
If $(\mathcal{A}, \bullet)$ is a trivial algebra, then it is immediate that $(\mathcal{A}, \bullet)$ is associative.
If $(\mathcal{A}, \bullet)$ is non-trivial, then by Corollary~\ref{co1}, the product $\bullet$ is associative.
\end{proof}

\begin{Proposition}\label{isometric-classification}
Every algebra given by the multiplication table and inner product in Corollary~\ref{co1} is a Euclidean left-alternative algebra. Moreover, two such algebras are isometrically isomorphic if and only if $s=s'$, $n=n'$, and there exists a permutation $\sigma$ of $\{1,\ldots,n\}$ such that
\[
(\lambda_i,t_i)=(\lambda'_{\sigma(i)},t'_{\sigma(i)})\qquad(1\leq i\leq n),
\]
where $t_i=0$ for $i>p$, and similarly for the primed parameters. The trivial algebras are classified up to isometric isomorphism by their dimension.
\end{Proposition}

\begin{proof}
Write $B_0=\operatorname{span}\{\mu_1,\ldots,\mu_s\}$ and $B_i=\operatorname{span}\{f_{i1},\ldots,f_{it_i}\}$, with $B_i=\{0\}$ if $i>p$. For elements
\[
x=b_0+\sum_{i=1}^n(\alpha_i e_i+b_i),\qquad
y=c_0+\sum_{i=1}^n(\beta_i e_i+c_i),
\]
the prescribed multiplication is
\[
x\bullet y=\sum_{i=1}^n\alpha_i(\beta_i e_i+c_i).
\]
For $z=d_0+\sum_i(\gamma_i e_i+d_i)$, both $(x\bullet y)\bullet z$ and $x\bullet(y\bullet z)$ equal $\sum_i\alpha_i\beta_i(\gamma_i e_i+d_i)$. Thus the algebra is associative. Moreover, $\Ll_x$ acts as the scalar $\alpha_i$ on $\mathbb R e_i\oplus B_i$ and vanishes on $B_0$. These subspaces are mutually orthogonal, so $\Ll_x$ is symmetric. This proves the converse.

Now let $F:\A\to\A'$ be an isometric algebra isomorphism. It preserves the left and two-sided annihilators and sends $N_\ell(\A)^\perp$ onto $N_\ell(\A')^\perp$. In particular, $s=s'$ and $n=n'$. The algebra $N_\ell(\A)^\perp$ is isomorphic to the coordinatewise algebra $\mathbb R^n$, whose primitive idempotents are precisely $e_1,\ldots,e_n$: an idempotent has coordinates in $\{0,1\}$ and is primitive exactly when one coordinate equals $1$. Hence $F(e_i)=e'_{\sigma(i)}$ for some permutation $\sigma$. Isometry gives $\lambda_i=\lambda'_{\sigma(i)}$. Since
\[
B_i=\{b\in N_\ell(\A)\mid e_i\bullet b=b\},
\]
we have $F(B_i)=B'_{\sigma(i)}$, and consequently $t_i=t'_{\sigma(i)}$. Conversely, given such a permutation, choose isometries $B_i\to B'_{\sigma(i)}$ and $B_0\to B'_0$, and extend them by $e_i\mapsto e'_{\sigma(i)}$. The resulting linear map is an isometry and preserves the multiplication by the displayed product formula. Finally, any two real inner product spaces of the same dimension are isometric, and for trivial algebras such an isometry automatically preserves the product.
\end{proof}

\begin{Corollary}
Let $(\mathcal{A}, \bullet, \prs)$ be a non-trivial commutative associative algebra endowed with a positive-definite associative bilinear form. Then, there exists a basis of $\mathcal{A}$ given by:
\[
\{ e_1, \dots, e_n,\,\mu_1, \dots, \mu_s\}
\]
where $n \in \mathbb{N}^*$ and $s \in \mathbb{N}$ (with the convention $s=0$ meaning that the family of vectors $\{\mu_1,\cdots, \mu_s\}$ is empty). With respect to this basis, the only nonzero products $\bullet$ and the only nonzero values of the bilinear form $\prs$ occur among the following:
\[
e_i \bullet e_i = e_i, \quad
\langle e_i, e_i \rangle = \lambda_i, \quad
\langle \mu_r, \mu_r \rangle = 1,
\]
for all $1 \leq i \leq n$ and $1 \leq r \leq s$, where $\lambda_i > 0$.
\end{Corollary}

\begin{proof}
By commutativity and associativity of the bilinear form, each $\Ll_u$ is symmetric. Hence $(\A,\bullet,\prs)$ is a Euclidean left-alternative algebra. In the normal form of Corollary~\ref{co1}, commutativity forces each $B_i$ to vanish, since $e_i\bullet b=b$ but $b\bullet e_i=0$ for $b\in B_i$. The result follows.
\end{proof}

\begin{Proposition}
Let $(\mathcal{A},\bullet,\prs)$ be a non-commutative Euclidean left-alternative algebra. Then, $\mathcal{A}^- := (\mathcal{A}, [ \cdot, \cdot ])$ is a two-step solvable Lie algebra but is not nilpotent. Moreover, there exists a basis:
\[
\{\mu_1, \dots, \mu_s, \; e_1, \dots, e_n, \; f_{11}, \dots, f_{1t_1}, \; \dots, \; f_{p1}, \dots, f_{pt_p}\}
\]
of $\mathcal{A}$, where $n,t_1,\ldots,t_p\in\mathbb{N}^*$, $1\leq p\leq n$, and $s \in \mathbb{N}$ (with the convention that $\{\mu_1, \dots, \mu_s\}$ is empty if $s = 0$), such that the only nonzero Lie brackets (up to skew-symmetry), the only nonzero Levi-Civita products $\star$, and the only nonzero values of the bilinear form $\prs$ are given by:
\begin{equation}\label{crochet}
\begin{aligned}
{}[e_k,f_{kt}]&=f_{kt},& f_{kt}\star f_{kt}&=\frac{1}{\lambda_k}e_k,& f_{kt}\star e_k&=-f_{kt},\\
\langle e_i,e_i\rangle&=\lambda_i,&\langle f_{kt},f_{kt}\rangle&=1,&\langle\mu_r,\mu_r\rangle&=1.
\end{aligned}
\end{equation}
for all $1 \leq i \leq n$, $1 \leq k \leq p$, $1 \leq t \leq t_k$, and $1 \leq r \leq s$, with $\lambda_i > 0$.  

In particular, the Riemannian Lie algebra $(\mathcal{A}^-, \prs)$ is not flat.
\end{Proposition}

\begin{proof}
By Corollary~\ref{associative}, $(\mathcal{A}, \bullet)$ is associative, hence $\mathcal{A}^- := (\mathcal{A}, [ \cdot, \cdot ])$ is a Lie algebra.  
From Corollary~\ref{co1}, the only nonzero brackets, up to skew-symmetry, are $[e_k,f_{kt}]=f_{kt}$. Consequently, the derived algebra is the nonzero abelian subspace spanned by the $f_{kt}$, and $\A^-$ has derived length exactly two. Furthermore, $[\A^-,[\A^-,\A^-]]=[\A^-,\A^-]$, so $\A^-$ is not nilpotent. By the Koszul formula \eqref{LV},
\[
2\langle f_{kt}\star f_{kt},e_k\rangle
=-\langle[f_{kt},e_k],f_{kt}\rangle+\langle[e_k,f_{kt}],f_{kt}\rangle=2.
\]
Orthogonality and the remaining instances of \eqref{LV} give $f_{kt}\star f_{kt}=\lambda_k^{-1}e_k$ and $f_{kt}\star e_k=-f_{kt}$; all other products of basis vectors vanish.  

Finally, since $\ell_{e_k}=0$ and $[f_{kt},e_k]=-f_{kt}$, the curvature convention fixed in the introduction gives:
\[
\mathfrak{R}(f_{kt},e_k)f_{kt}=-\ell_{f_{kt}}(f_{kt})=-f_{kt}\star f_{kt}=-\frac{1}{\lambda_k}e_k\neq0,
\]
which shows that $(\mathcal{A}^-, \prs)$ is not flat.
\end{proof}

\begin{Remark}
Recall that in \cite{Milnor}, J. Milnor showed that every flat Riemannian Lie algebra has derived length at most two. The above proposition gives examples of non-flat two-step solvable Euclidean Lie algebras.
\end{Remark}

We conclude with a geometric application of the algebraic classification. The relevant constant-curvature Riemannian Lie groups are the classical examples described in \cite[Special Example 1.7]{Milnor}. The contribution here is to identify the factors and their metric parameters directly from the structure of the original left-alternative algebra.

\begin{Theorem}\label{hyperbolic-product}
Let $(\A,\bullet,\prs)$ be a Euclidean left-alternative algebra, and let $G$ be the simply connected Lie group with Lie algebra $\A^-$. Equip $G$ with the left-invariant Riemannian metric induced by $\prs$. With the notation of Corollary~\ref{co1}, and allowing $p=0$ in the commutative case, $G$ is isometric, as a Riemannian manifold, to the Riemannian product
\[
\mathbb R^{s+n-p}\times\prod_{i=1}^{p}\mathbb H^{t_i+1}_{-1/\lambda_i},
\]
where $\mathbb H^d_K$ denotes the simply connected complete real hyperbolic space of dimension $d$ and constant sectional curvature $K<0$. For a trivial algebra, take $n=p=0$ and $s=\dim\A$. In particular, the metric is flat if and only if $(\A,\bullet)$ is commutative.
\end{Theorem}

\begin{proof}
By the classification, the Lie algebra $\A^-$ is the orthogonal direct sum of the central abelian ideal $\mathfrak z_0=\operatorname{span}\{\mu_1,\ldots,\mu_s,e_{p+1},\ldots,e_n\}$ and the ideals $\mathfrak g_i=\mathbb R e_i\oplus B_i$, for $1\leq i\leq p$, where $B_i=\operatorname{span}\{f_{i1},\ldots,f_{it_i}\}$. The only nonzero brackets in $\mathfrak g_i$, up to skew-symmetry, are $[e_i,b]=b$ for $b\in B_i$. Let $G_0$ and $G_i$ be the simply connected Lie groups with Lie algebras $\mathfrak z_0$ and $\mathfrak g_i$, respectively. By the correspondence between simply connected Lie groups and Lie algebras (see \cite[Chapter 20]{LeeSmooth}), the direct sum decomposition into ideals integrates to a Lie group isomorphism $G_0\times G_1\times\cdots\times G_p\simeq G$. Orthogonality at the identity and left invariance imply that its metric is the product metric. The factor corresponding to $\mathfrak z_0$ is Euclidean of dimension $s+n-p$.

Fix $i$ and put $a=e_i/\sqrt{\lambda_i}$ and $c=1/\sqrt{\lambda_i}$. Then $a$ is a unit vector orthogonal to $B_i$, and $[a,b]=cb$ for all $b\in B_i$. Define the linear form $\varphi_i$ on $\mathfrak g_i$ by $\varphi_i(x)=c\langle a,x\rangle$. For all $x,y\in\mathfrak g_i$, the bracket is
\[
[x,y]=\varphi_i(x)y-\varphi_i(y)x.
\]
Hence, by \cite[Special Example 1.7, pp.~298--299]{Milnor}, the left-invariant metric on $G_i$ has constant sectional curvature $-\|\varphi_i\|^2=-c^2=-1/\lambda_i$. For completeness, we also give the direct calculation in our curvature convention. By the Koszul formula~\eqref{LV}, for all $b,d\in B_i$,
\[
a\star a=a\star b=0,\qquad b\star a=-cb,\qquad b\star d=c\langle b,d\rangle a.
\]
In particular, $\ell_a=0$. Using the curvature convention fixed in the introduction, we obtain, for all $b,d,h\in B_i$,
\begin{align*}
\mathfrak R(a,b)a&=-c^2b,&
\mathfrak R(a,b)d&=c^2\langle b,d\rangle a,\\
\mathfrak R(b,d)a&=0,&
\mathfrak R(b,d)h&=c^2(\langle d,h\rangle b-\langle b,h\rangle d).
\end{align*}
Indeed, $\mathfrak R(a,b)=c\ell_b$ and $\mathfrak R(b,d)=-[\ell_b,\ell_d]$, which give the displayed identities upon substitution of the Levi-Civita products. By trilinearity and skew-symmetry in the first two arguments, for all $x,y,z\in\mathfrak g_i$ we have
\[
\mathfrak R(x,y)z=c^2(\langle y,z\rangle x-\langle x,z\rangle y).
\]
With our curvature convention, the sectional curvature of the plane spanned by linearly independent $x,y$ is
\[
K(x,y)=\frac{\langle\mathfrak R(x,y)x,y\rangle}{\langle x,x\rangle\langle y,y\rangle-\langle x,y\rangle^2}=-c^2=-\frac1{\lambda_i}.
\]
The Riemannian Lie group $G_i$ is geodesically complete: local existence of geodesics and compactness of the unit sphere at the identity give a uniform positive existence time for unit-speed geodesics starting there, and left translations transport this uniform time to every point. Thus unit-speed geodesics can be continued indefinitely in both directions. Since $G_i$ is simply connected, the Killing--Hopf theorem (see \cite[Chapter 12]{LeeRiemannian}) implies that $G_i$, viewed as a Riemannian manifold, is isometric to $\mathbb H^{t_i+1}_{-1/\lambda_i}$. This proves the asserted product decomposition. Finally, the metric is flat precisely when $p=0$, which, by the multiplication table, is equivalent to commutativity of $(\A,\bullet)$.
\end{proof}

\section*{Data Availability}
Data sharing is not applicable to this article, as no datasets were
generated or analyzed during the current study.

\section*{Author Contributions}
All authors contributed equally to this work and approved the final
manuscript.

\section*{Competing Interests}
The authors declare that they have no competing interests.

\section*{Funding}
H. E. O. was supported by grant NYUAD-065.

\bibliographystyle{elsarticle-num}

\end{document}